# Bimonoidal Categories, $E_n$-Monoidal Categories, and Algebraic $K$-Theory

Volumes I and II by Donald Yau

Volume III by Niles Johnson and Donald Yau

This is an abridged version of the book for arxiv.org that contains the full list of open questions together with front and back matters.
The complete pdf version is available on the authors' web pages:

https://nilesjohnson.net/En-monoidal
https://u.osu.edu/yau.22/main

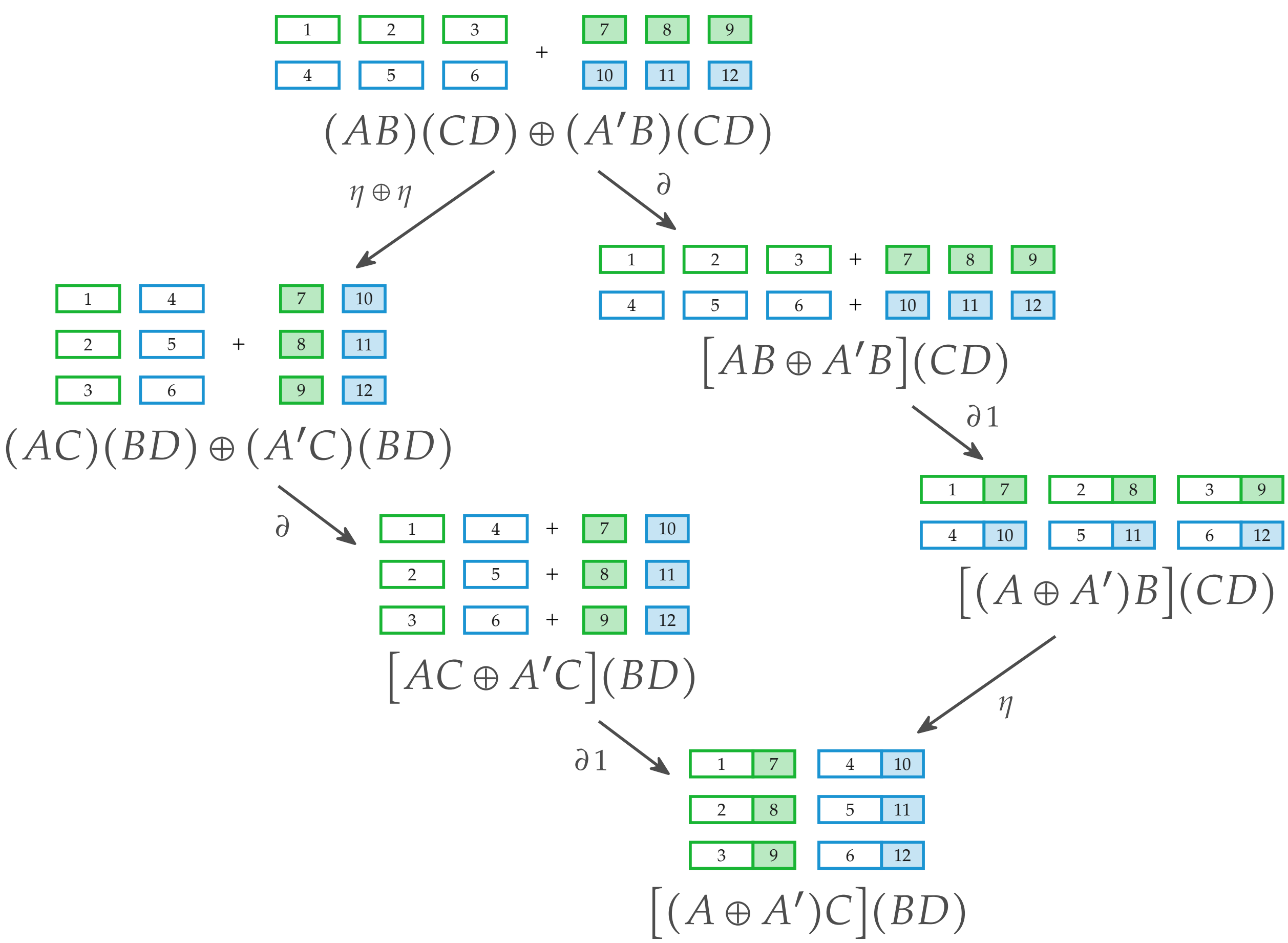

Department of Mathematics, The Ohio State University at Newark,
1179 University Drive, Newark, OH 43055, USA
*Email address*: johnson.5320@osu.edu, yau.22@osu.edu

ABSTRACT. Bimonoidal categories are categorical analogues of rings without additive inverses. They have been actively studied in category theory, homotopy theory, and algebraic $K$-theory since around 1970. There is an abundance of new applications and questions of bimonoidal categories in mathematics and other sciences. This work provides a unified treatment of bimonoidal and higher ring-like categories, their connection with algebraic $K$-theory and homotopy theory, and applications to quantum groups and topological quantum computation. With ample background material, extensive coverage, detailed presentation of both well-known and new theorems, and a list of open questions, this work is a user friendly resource for beginners and experts alike.

Part I.1 proves in detail Laplaza's two coherence theorems and May's strictification theorem of symmetric bimonoidal categories, as well as their bimonoidal analogues. Part I.2 proves Baez's Conjecture on the existence of a bi-initial object in a 2-category of symmetric bimonoidal categories. The next main theorem states that a matrix construction, involving the matrix product and the matrix tensor product, sends a symmetric bimonoidal category with invertible distributivity morphisms to a symmetric monoidal bicategory, with no strict structure morphisms in general.

Part II.1 studies braided bimonoidal categories, with applications to quantum groups and topological quantum computation. It is proved that the categories of modules over a braided bialgebra, of Fibonacci anyons, and of Ising anyons form braided bimonoidal categories. Two coherence theorems for braided bimonoidal categories are proved, confirming the Blass-Gurevich conjecture. The rest of this part discusses braided analogues of Baez's Conjecture and the monoidal bicategorical matrix construction in Part I.2. Part II.2 studies ring and bipermutative categories in the sense of Elmendorf-Mandell, braided ring categories, and $E_n$-monoidal categories, which combine $n$-fold monoidal categories with ring categories.

Part III.1 is a detailed study of enriched monoidal categories, pointed diagram categories, and enriched multicategories. Using the machinery in Part III.1, Part III.2 discusses the rich interconnection between the higher ring-like categories in Part II.2, homotopy theory, and algebraic $K$-theory. Starting with a chapter on homotopy theory background, the first half of this part constructs the Segal $K$-theory functor and the Elmendorf-Mandell $K$-theory multifunctor from permutative categories to symmetric spectra. The second half applies the $K$-theory multifunctor to small ring, bipermutative, braided ring, and $E_n$-monoidal categories to obtain, respectively, strict ring, $E_\infty$-, $E_2$-, and $E_n$-symmetric spectra. Appendix III.A discusses open questions related to the topics of this work.

# Contents

*Reference Convention*. There are six parts in this work, with two parts in each volume. In each of Volumes II and III, the numbering of parts, chapters, and pages restarts at 1. Page numbers are prefixed with a volume number so that, for example, page $n$ of Volume II is printed II.$n$. References (chapters, definitions, theorems, equations, etc.) to a different volume are preceded by I, II, or III according to the volume. For example, $a.b.c$ in Volume I is denoted by I.$a.b.c$ in Volumes II and III, and Chapter $N$ in Volume III is denoted by Chapter III.$N$ in Volumes I and II.

# Preface

Bimonoidal categories are categorical analogues of rings without additive inverses. They have been actively studied in category theory, homotopy theory, and algebraic $K$-theory since around 1970. There is an abundance of new applications and questions of bimonoidal categories in mathematics and other sciences. This work provides the first unified treatment of bimonoidal and higher ring-like categories, their connection with algebraic $K$-theory and homotopy theory, and applications to quantum groups and topological quantum computation. With ample background material, extensive coverage, detailed presentation of both well-known and new theorems, and a list of open questions, this work is a user friendly resource for beginners and experts alike.

### Bimonoidal and $E_n$-Monoidal Categories

A *bimonoidal category* $\mathsf{C}$ is a categorical analogue of a rig, which is a ring without additive inverses. In place of the rig addition is a symmetric monoidal structure

$$\big(\mathsf{C}, \oplus, 0, \alpha^{\oplus}, \lambda^{\oplus}, \rho^{\oplus}, \xi^{\oplus}\big),$$

called the additive structure. In place of the rig multiplication is a monoidal structure

$$\big(\mathsf{C}, \otimes, \mathbb{1}, \alpha^{\otimes}, \lambda^{\otimes}, \rho^{\otimes}\big),$$

called the multiplicative structure. The multiplicative zero property in a rig, $0x = 0 = x0$, is replaced by two natural isomorphisms

$$0 \otimes A \xrightarrow[\cong]{\lambda^{\bullet}_A} 0 \xleftarrow[\cong]{\rho^{\bullet}_A} A \otimes 0$$

for objects $A \in \mathsf{C}$, called the left and right multiplicative zeros. The categorical analogues of the distributivity properties in a rig,

$$x(y+z) = xy + xz \quad \text{and} \quad (x+y)z = xz + yz,$$

are two natural monomorphisms

$$A \otimes (B \oplus C) \xrightarrow{\delta^l_{A,B,C}} (A \otimes B) \oplus (A \otimes C)$$

$$(A \oplus B) \otimes C \xrightarrow{\delta^r_{A,B,C}} (A \otimes C) \oplus (B \otimes C)$$

called the left and right distributivity morphisms. A *symmetric (braided) bimonoidal category* has the additional structure of a natural isomorphism

$$A \otimes B \xrightarrow[\cong]{\xi^{\otimes}_{A,B}} B \otimes A$$

that makes the multiplicative structure into a symmetric (braided) monoidal category. These data are required to satisfy a finite list of axioms that (i) are checkable in practice and (ii) ensure that (symmetric/braided) bimonoidal categories have good coherence and other categorical properties. If the distributivity morphisms $\delta^l$ and $\delta^r$ are natural isomorphisms, then we call $\mathsf{C}$ *tight*.

For example, the category $\mathsf{Vect}^{\mathbb{C}}$ of finite dimensional complex vector spaces is a tight symmetric bimonoidal category with the usual direct sum and tensor product. More generally, each distributive symmetric monoidal category is a tight symmetric bimonoidal category. The nonnegative integers and permutations form the objects and the morphisms of a tight symmetric bimonoidal category $\Sigma$, called the finite ordinal category. May's bipermutative categories, with the additional axiom $\xi^{\otimes}_{-,0} = \mathrm{Id}$, are tight symmetric bimonoidal categories.

We introduce higher analogues of bimonoidal categories called *$E_n$-monoidal categories* for $n \geq 1$. On top of a symmetric monoidal structure $\oplus$, an $E_n$-monoidal category has interacting monoidal structures $\{\otimes_i\}_{i=1}^n$ and factorization morphisms, along with appropriate compatibility axioms. The notion of an $E_n$-monoidal category simultaneously generalizes

- ring categories in the sense of Elmendorf-Mandell **[EM06]** with $n = 1$,
- bipermutative categories **[EM06]** with $n \geq 2$,
- braided bimonoidal categories in the sense of Richter **[Ric10]** with $n = 2$, and
- $n$-fold monoidal categories as in **[BFSV03]** with the monoidal structures $\{\otimes_i\}_{i=1}^n$.

The braided bimonoidal categories in this work, which are studied in Part II.1, are more general than Richter's; see Note II.9.7.5.

## Applications

Bimonoidal and higher ring-like categories are at the crossroad of category theory, algebraic $K$-theory, homotopy theory, and other sciences.

*Category Theory*. As categorifications of rigs, bimonoidal categories and their symmetric and braided analogues are interesting categorical objects in their own right. Laplaza **[Lap72a, Lap72b]** defined symmetric bimonoidal categories and proved two coherence theorems in the form of commutative formal diagrams.

- Laplaza's First Coherence Theorem I.3.9.1 is analogous to Mac Lane's Coherence Theorem I.1.3.3 for monoidal categories. This theorem is widely cited and used in the literature, and has far-reaching consequences, including the Strictification Theorems I.5.4.6 and I.5.4.7 for tight symmetric bimonoidal categories.
- Laplaza's Second Coherence Theorem I.4.4.3 is analogous to the Joyal-Street Coherence Theorem II.1.6.3 for braided monoidal categories. We prove a braided analogue in Theorem II.5.4.4, which confirms the Blass-Gurevich Conjecture **[BG20a]**.

Ring and bipermutative categories in the sense of Elmendorf-Mandell are different categorifications of rigs and commutative rigs, with factorization morphisms instead of distributivity morphisms. These categories and their braided and higher analogues are discussed in Part II.2.

*Algebraic K-Theory and Homotopy Theory*. Ring-like categories are the inputs of algebraic $K$-theory functors due to

- Segal [**Seg74**],
- May [**May77, May78, May82, May09a**], and
- Elmendorf-Mandell [**EM06, EM09**],

among others, that produce structured ring spectra, which are among the most important objects in homotopy theory. The deep connection between category theory, algebraic $K$-theory, and homotopy theory is the subject of Part III.2.

*Other Sciences*. Due to the ubiquity of ring-like structures and categories, bimonoidal categories are increasingly applied in other sciences, including computer science. See Note I.2.7.5 for related references. Section I.2.6 contains an application of symmetric bimonoidal categories to reversible programming of finite types. Applications of braided bimonoidal categories to quantum group theory and topological quantum computation are discussed in Chapter II.3.

## Purpose

This work is a systematic treatise of the following aspects of bimonoidal categories and their symmetric, braided, and higher analogues:

- 1- and 2-dimensional category theory,
- applications of braided bimonoidal categories to quantum groups and topological quantum computation, and
- algebraic $K$-theory.

Here are some of the highlights.

**Category Theory:**

*Bimonoidal Coherence*. Coherence and strictification theorems for symmetric bimonoidal categories have been extensively used in the literature since the 1970s. In addition to providing completely detailed proofs of these theorems, we also correct some subtle and nontrivial inaccuracies in the original proofs and statements in the coherence theorems in [**Lap72a, Lap72b**] that have never been made explicit before. See Sections I.3.11 and I.4.7 for related discussion. These theorems and their bimonoidal analogues are discussed in detail in Part I.1.

*The Blass-Gurevich Conjecture*. For braided bimonoidal categories, our coherence and strictification theorems are new and confirm a conjecture of Blass and Gurevich [**BG20a**]. Braided bimonoidal categories are discussed in Part II.1. Their coherence theorems are discussed in Chapters II.5 and II.6.

*Centers*. Monoidal, braided monoidal, and symmetric monoidal categories are connected by the Drinfeld center and the symmetric center. Bimonoidal and ring categorical analogues of these center constructions are discussed in Chapters II.4 and II.9.

*$E_n$-Monoidal Categories*. Just as $n$-fold monoidal categories contain strict (braided) monoidal and permutative categories, in Part II.2 we show that our $E_n$-monoidal categories contain ring and bipermutative categories in [**EM06**] and braided bimonoidal categories in [**Ric10**].

**Two-Dimensional Categories:** Parts I.2 and II.1 contain several new theorems on the close connection between bimonoidal categories and 2-dimensional categories.

*Baez's Conjecture*. In Chapter I.7, we prove a conjecture of Baez [**Bae18**] on the existence of a bi-initial object in a 2-category of symmetric bimonoidal categories. Section I.7.9 discusses the relationship between our version of Baez's Conjecture for symmetric bimonoidal categories and a more restricted version in [**CDH∞, Elg21**] for rig categories.

*Monoidal Bicategorification*. In Chapter I.8, we show that each tight symmetric bimonoidal category $\mathsf{C}$ yields a symmetric monoidal bicategory $\mathsf{Mat}^{\mathsf{C}}$ whose 1-/2-cells are matrices in the objects/morphisms in $\mathsf{C}$. Its other structures involve the matrix product and the matrix tensor product, and they are not strict in general. Coordinatized 2-vector spaces form such a symmetric monoidal bicategory. Chapter I.9 discusses the analogues in which $\mathsf{Mat}^{\mathsf{C}}$ is a permutative Gray monoid or a permutative 2-category.

*Braided Versions*. In Chapters II.7 and II.8, we prove braided analogues of, respectively, Baez's Conjecture (Chapter I.7) and the monoidal bicategorification theorem (Chapter I.8) for braided bimonoidal categories.

**Applications:**

*Quantum Groups*. The first part of Chapter II.3 extends a well-known fact in quantum group theory. We observe that the category of left modules over a braided bialgebra, which is also known as a quasitriangular bialgebra in the literature, is a tight braided bimonoidal category.

*TQC*. The second part of Chapter II.3 discusses applications of braided bimonoidal categories to topological quantum computation (TQC). We prove that the Fibonacci anyons and the Ising anyons, which are two of the most important models in TQC, are both tight braided bimonoidal categories.

*Programming*. Section I.2.6 is a brief illustration that symmetric bimonoidal categories naturally arise in reversible programming of finite types.

**Algebraic *K*-Theory:** Part III.2 discusses the interplay between $E_n$-monoidal categories (Part II.2), homotopy theory, and algebraic $K$-theory. The first half of Part III.2 discusses in detail

- the Segal $K$-theory functor and
- the Elmendorf-Mandell $K$-theory multifunctor

from small permutative categories to symmetric spectra. Our treatment corrects a subtle but nontrivial inaccuracy in [**EM09**, Theorem 1.3] and some other statements about expanding the domain of the $K$-theory multifunctor. See Note III.10.8.2 for a complete explanation.

The second half of Part III.2 applies the $K$-theory multifunctor to the $E_n$-monoidal categories in Part II.2 to produce structured ring spectra. We prove in detail that the $K$-theory multifunctor sends

- small ring categories to strict ring symmetric spectra,
- small bipermutative categories to $E_\infty$-symmetric spectra,
- small braided ring categories to $E_2$-symmetric spectra, and

- small $E_n$-monoidal categories to $E_n$-symmetric spectra for $2 \leq n < \infty$.

The strict ring and $E_\infty$ cases are from **[EM06, EM09]**. The $1 < n < \infty$ cases are new results. These theorems and the detailed construction of the algebraic $K$-theory (multi)functors use a substantial amount of enriched monoidal category and multicategory theory, which is the subject of Part III.1.

**Audience and Features**

This work is aimed at graduate students and advanced researchers with an interest in category theory, homotopy theory, algebraic $K$-theory, and their applications. Below are some features that make this work a unique and user friendly resource.

**Unified Presentation:** The literature on bimonoidal categories, higher ring-like categories, enriched monoidal categories, multicategories, and their connection with algebraic $K$-theory, homotopy theory, and the sciences is scattered across many journal articles over several decades, with varying definitions, notations, and terminology. This work presents these topics in a unified manner, with both well-known and new theorems.

**Background Material:** To make this work self-contained and to bring the reader quickly up to speed, there is extensive background material on

- basic category theory (Chapter I.1),
- 2-dimensional categories (Chapter I.6),
- braided structures (Chapter II.1),
- abelian categories (Section II.2.3),
- braided, also known as quasitriangular, bialgebras (Section II.3.1),
- enriched monoidal categories (Chapters III.1, III.2, and III.3),
- pointed objects and pointed diagram categories (Chapter III.4),
- enriched multicategories (Chapters III.5 and III.6), and
- homotopy theory (Chapter III.7).

These chapters and sections form a substantial portion of this work.

**Open Questions:** Appendix III.A discusses open questions related to the topics of this work. The reader is encouraged to take advantage of these open questions and use them as a springboard to read the main text.

**Detailed Discussion:** This work contains many highly detailed and carefully structured proofs for both known and new theorems. For each major result, our discussion has much more detail than one would normally find in the literature. Our detailed discussion has several purposes.

*Exercises with Solutions*. Our detailed presentation makes the material accessible to a diverse audience, including those who are new to bimonoidal and higher ring-like categories and algebraic $K$-theory. Students are encouraged to regard the numerous detailed proofs as exercises with full solutions. Each result whose proof has many different parts has been carefully structured to make it easy for the reader to jump forward and backward.

*Axioms*. Symmetric bimonoidal categories are defined by 24 axioms, and the list of axioms for (braided) bimonoidal categories is similarly substantial. Our detailed discussion helps the reader see exactly where these axioms are used and why they are needed.

*Laplaza's Theorems*. The Coherence Theorems I.3.9.1 and I.4.4.3 for symmetric bimonoidal categories are central results in this subject that have been cited and used numerous times in the literature. Their original proofs given by Laplaza in [**Lap72a, Lap72b**] were written in outline form, with much detail and some cases in the proofs completely omitted. Moreover, Laplaza's original proofs and statements of these theorems have several subtle and nontrivial inaccuracies that have never been made explicit before and are not easy to spot. For both archival and educational purposes, we present fully detailed proofs of these theorems and correct the inaccuracies. Sections I.3.11 and I.4.7 have more related discussion.

*K-Theory Multifunctors*. The $K$-theory multifunctors in Chapters III.9 and III.10, due to Elmendorf-Mandell [**EM06, EM09**], are fundamental constructions for multiplicative structure of algebraic $K$-theory spectra. They are essential for our development of $E_n$-monoidal symmetric spectra from corresponding structure on small permutative categories. We use the theory of enriched monoidal categories and enriched multicategories from Part III.1 to give complete explanations of the constructions and their properties. This treatment corrects an inaccuracy in the statement of [**EM09**, Theorem 1.3] and some other statements about expanding the domain of the $K$-theory multifunctor. The basic issue has to do with monoidal units and, to the authors' knowledge, has not been previously explained. See Note III.10.8.2 for further discussion.

**Reading Guide:** In addition to a detailed introduction, almost every chapter has a brief *Reading Guide* that provides an alternative to reading that chapter linearly. Our presentation in the main text follows a straightly logical order and has a lot of detail. By following the reading guide, it is possible to first obtain a bird's-eye view of that chapter before digesting all the detail.

**Motivation and Explanation:** Main definitions and results are often preceded by discussion that motivates the upcoming definitions and proofs. Whenever useful, definitions and results are followed by a detailed explanation that interprets and unpacks the various components. In the text, these are clearly marked as *Motivation* and *Explanation*, respectively. Examples include Motivation I.2.1.1, Explanation I.2.4.7, and Section I.4.1.

**Organization:** There are extensive cross-references throughout the text. In addition to a detailed index, there are lists of main facts and notations, each organized by chapters. While the text follows a strictly logical order, it is not necessary to read the chapters in a linear order. The reader can jump straight to a section and use the extensive cross-references to fill in the necessary definitions and facts.

### Related Literature

Here we list a selection of general references for background or further reading. The Notes section at the end of each chapter provides additional references for the content of that chapter.

Basic category theory**:** [**Awo10, BK00, Gra18, Lei14, Rie16, Rom17, Sim11**].
More advanced category theory**:** [**Bor94a, Bor94b, Bor94c, ML98, Mit65, Sch72**].

Abelian categories**:** [**EGNO15, Fre03, Mit65**].
Enriched categories**:** [**Bor94b, Cru09, For04, Kel05**].
Ends and coends**:** [**Day70, DK69, Lor21**].
2-dimensional categories**:** We highly recommend [**JY21**].
Multicategories**:** [**Fre17, MSS02, May72, Yau16**].
Homotopy theory**:** [**BR20, May99, MP12, Mil20, Ric20, Rie14**].
Simplicial homotopy theory**:** [**Cur71, GZ67, GJ09, May92**].
Algebraic $K$-theory**:** [**Mil71, Qui73, Ros95, Wal85, Wei13**].

## Part and Chapter Summaries

### Part I.1: Symmetric Bimonoidal Categories

This part studies symmetric bimonoidal categories and bimonoidal categories (Chapter I.2). It presents highly detailed proofs of Laplaza's Coherence Theorems for symmetric bimonoidal categories (Chapters I.3 and I.4), May's Strictification Theorem for tight symmetric bimonoidal categories (Chapter I.5), and their non-symmetric analogues for bimonoidal categories. The only prerequisite for this part is some basic knowledge of category theory, which is summarized in Chapter I.1.

### Part I.2: Bicategorical Aspects of Symmetric Bimonoidal Categories

Applying Laplaza's Coherence Theorems, this part proves several new theorems on the connection between symmetric bimonoidal categories and bicategories. All the necessary definitions of 2-dimensional category theory are summarized in Chapter I.6. The first main result is a confirmation of Baez's Conjecture (Chapter I.7) that proves the existence of a bi-initial object in a 2-category of symmetric bimonoidal categories. Chapter I.8 proves that a matrix construction $\mathsf{Mat}^{\mathsf{C}}$ sends each tight symmetric bimonoidal category to a symmetric monoidal bicategory. With further strictness assumptions, $\mathsf{Mat}^{\mathsf{C}}$ is a permutative Gray monoid or even a permutative 2-category (Chapter I.9).

### Part II.1: Braided Bimonoidal Categories

Starting with a preliminary chapter on the braid groups and braided monoidal categories, this part is a detailed study of braided bimonoidal categories (Chapter II.2), which are strictly more general than Richter's [**Ric10**] and the BD categories of Blass-Gurevich [**BG20a**]. This part discusses applications to quantum groups and topological quantum computation (Chapter II.3), bimonoidal centers (Chapter II.4), coherence and strictification of braided bimonoidal categories (Chapters II.5 and II.6), and the braided versions of Baez's Conjecture and the matrix construction (Chapters II.7 and II.8). Our coherence and strictification theorems confirm the Blass-Gurevich Conjecture. The main theorems in Parts I.1 and I.2 are used in this part.

### Part II.2: $E_n$-Monoidal Categories

This part studies a closely related variant of bimonoidal categories, called ring categories, and their bipermutative, braided, and higher analogues, called $E_n$-monoidal categories. Ring and bipermutative categories are due to Elmendorf-Mandell [**EM06, EM09**]. An $E_n$-monoidal category combines $n$ ring categories with a common additive structure and an $n$-fold monoidal category as in [**BFSV03**]. The categories in this part are applied in Part III.2 to obtain $E_n$-symmetric spectra via algebraic $K$-theory. This part is independent of the earlier parts, except for some definitions and statements of theorems.

Part III.1: Enriched Monoidal Categories and Multicategories

To prepare for Part III.2, this part lays the groundwork on enriched monoidal categories (Chapters III.1, III.2, and III.3), smash products (Chapter III.4), and multicategories (Chapters III.5 and III.6). In addition to their roles in the Segal $K$-theory functor and the Elmendorf-Mandell $K$-theory multifunctor, the detailed discussion of enriched monoidal categories—including change of enrichment, coherence, self-enrichment, and the Enriched Yoneda Lemma—and multicategories is also of independent interest. These chapters assume only a basic knowledge of monoidal categories, as summarized in Section III.1.1.

Part III.2: Algebraic $K$-Theory

This part studies the interconnection between $E_n$-monoidal categories (Part II.2), homotopy theory (Chapter III.7), and algebraic $K$-theory. The first half discusses in detail the Segal $K$-theory functor (Chapter III.8) and the Elmendorf-Mandell $K$-theory multifunctor (Chapters III.9 and III.10) from small permutative categories to symmetric spectra. The second half (Chapters III.11, III.12, and III.13) applies the $K$-theory multifunctor to small ring, bipermutative, braided ring, and $E_n$-monoidal categories to obtain, respectively, strict ring, $E_\infty$-, $E_2$-, and $E_n$-symmetric spectra. These structured ring spectra are fundamental objects in homotopy theory. Our discussion shows how they arise from $E_n$-monoidal categories via algebraic $K$-theory.

In the main text, each chapter starts with a detailed introduction. A summary of each chapter follows.

**Part I.1: Symmetric Bimonoidal Categories**

Chapter I.1: Basic Category Theory

To make this book self-contained, this chapter reviews the basics of category theory, starting from the definitions of categories, functors, and natural transformations. Then it discusses adjunctions, equivalences of categories, (co)limits, (co)ends, and Kan extensions. The remaining sections review (symmetric) monoidal categories, (symmetric) monoidal functors, monoidal natural transformations, and their coherence theorems.

Chapter I.2: Symmetric Bimonoidal Categories

This chapter introduces symmetric bimonoidal categories and bimonoidal categories. Then we prove Laplaza's Theorem I.2.2.13 that says that half of the 24 symmetric bimonoidal category axioms are formal consequences of the other 12 axioms. The weaker bimonoidal analogue is Proposition I.2.2.14. The remaining sections discuss examples of symmetric bimonoidal categories, including distributive symmetric monoidal categories, the finite ordinal category $\Sigma$, a variant $\Sigma'$, and left and right bipermutative categories. The finite ordinal category $\Sigma$ is an important part of (i) the distortion category $\mathcal{D}$ (Chapter I.4) used in Laplaza's Second Coherence Theorem I.4.4.3, (ii) Baez's Conjecture (Chapter I.7), and (iii) the braided version of Baez's Conjecture (Chapter II.7). Section I.2.6 contains an application of symmetric bimonoidal categories to reversible programming of finite types.

Chapter I.3: Coherence of Symmetric Bimonoidal Categories

This chapter proves Laplaza's First Coherence Theorem I.3.9.1 for symmetric bimonoidal categories that satisfy a monomorphism assumption. This assumption is automatically satisfied if tightness—that is, the invertibility of the distributivity morphisms $\delta^l$ and $\delta^r$—is assumed, but the general form of this theorem only requires that the distributivity morphisms be natural monomorphisms. The analogue of this coherence theorem for bimonoidal categories is Theorem I.3.10.7. Section I.3.11 discusses the main differences between this chapter and Laplaza's original work in [**Lap72a**].

Chapter I.4: Coherence of Symmetric Bimonoidal Categories II

This chapter proves Laplaza's Second Coherence Theorem I.4.4.3 for symmetric bimonoidal categories that satisfy the same monomorphism assumption as in Theorem I.3.9.1. The analogue of this coherence theorem for bimonoidal categories is Theorem I.4.5.8. Section I.4.7 discusses the main differences between this chapter and Laplaza's original work in [**Lap72b**]. Both Coherence Theorems I.3.9.1 and I.4.4.3 say that some formal diagrams in certain symmetric bimonoidal categories commute. The first theorem has an assumption called *regularity* on the common domain of the two paths involved, which is analogous to Mac Lane's Coherence Theorem I.1.3.3 for monoidal categories. The second theorem has an assumption about the two paths themselves, which is reminiscent of the Joyal-Street Coherence Theorem II.1.6.3 for braided monoidal categories. In Chapter II.5, we observe that the second, but not the first, theorem has a braided analogue.

Chapter I.5: Strictification of Tight Symmetric Bimonoidal Categories

This chapter proves May's Strictification Theorem I.5.4.6 of *tight* symmetric bimonoidal categories to right bipermutative categories. The latter are tight symmetric bimonoidal categories whose additive structures and multiplicative structures are both permutative categories, and whose structure morphisms $\lambda^\bullet$, $\rho^\bullet$, $\delta^r$, and $\xi^{\otimes}_{-,0}$ are identities. Unlike the Coherence Theorems I.3.9.1 and I.4.4.3, the strictification theorem requires the tightness assumption. Our detailed proofs show exactly where the invertibility of $\delta^l$ and $\delta^r$ is used. Theorem I.5.4.7 is another version of the strictification theorem involving *left* bipermutative categories, in which $\delta^l$, instead of $\delta^r$, is the identity. Theorems I.5.5.11 and I.5.5.12 are the corresponding strictification results for tight bimonoidal categories. Section I.5.6 briefly discusses the history of related strictification theorems and claims.

**Part I.2: Bicategorical Aspects of Symmetric Bimonoidal Categories**

Chapter I.6: Definitions from Bicategory Theory

Without assuming any knowledge of 2-dimensional categories, in this chapter we review the basics of 2-/bicategories, pasting diagrams, lax functors, lax transformations, modifications, and adjunctions in bicategories. Then it reviews multiplicative structures, including monoidal bicategories, their braided, sylleptic, and symmetric analogues, the Gray tensor product for 2-categories, (permutative) Gray monoids, and permutative 2-categories. Most of these topics are discussed in detail in the book [**JY21**].

Chapter I.7: Baez's Conjecture

This chapter proves Baez's Conjecture (Theorems I.7.8.1 and I.7.8.3). Section I.7.1 defines a 2-category $\mathsf{Bi}_\mathsf{r}^{\mathsf{fsy}}$ with *flat* small symmetric bimonoidal categories as objects. Flatness (Definition I.3.9.9) is much weaker than tightness, and it guarantees that the Coherence Theorems I.3.9.1 and I.4.4.3 are applicable. The first version of Baez's Conjecture (Theorem I.7.8.1) says that the finite ordinal category $\Sigma$ is a lax bicolimit of the 2-functor $\varnothing \longrightarrow \mathsf{Bi}_\mathsf{r}^{\mathsf{fsy}}$. Another version is Theorem I.7.8.3, which says that the variant $\Sigma'$ of $\Sigma$ is also such a lax bicolimit. We emphasize that our proof of Baez's Conjecture does *not* use the Strictification Theorems I.5.4.6 and I.5.4.7. This allows us to use flat small symmetric bimonoidal categories in the 2-category $\mathsf{Bi}_{\mathsf{r}}^{\mathsf{fsy}}$, instead of the smaller class of tight ones. Section I.7.9 discusses the relationship between our version of Baez's Conjecture and the more restricted version in [**CDH∞, Elg21**] for rig categories, which are multiplicatively nonsymmetric and tight.

Chapter I.8: Symmetric Monoidal Bicategorification

This chapter proves Theorem I.8.15.4. It says that, for each tight symmetric bimonoidal category $\mathsf{C}$, a matrix construction $\mathsf{Mat}^{\mathsf{C}}$ is a symmetric monoidal bicategory, with no strict structure morphisms in general. Therefore, the construction $\mathsf{Mat}^{\mathsf{C}}$ is a direct connection between tight symmetric bimonoidal categories and symmetric monoidal bicategories. The objects in $\mathsf{Mat}^{\mathsf{C}}$ are nonnegative integers. Its 1-/2-cells are matrices whose entries are objects/morphisms in $\mathsf{C}$. The horizontal composition in the bicategory $\mathsf{Mat}^{\mathsf{C}}$ uses the usual matrix product. The monoidal composition in its monoidal bicategory structure uses the matrix tensor product, which is also known as the Kronecker product. The category of coordinatized 2-vector spaces, which is $\mathsf{Mat}^{\mathsf{C}}$ with $\mathsf{C} = \mathsf{Vect}^{\mathbb{C}}$, is such a symmetric monoidal bicategory. This chapter uses the Coherence Theorems I.3.9.1 and I.3.10.7 and the graph theoretic machinery in Chapter I.3, but neither the Coherence Theorem I.4.4.3 nor the Strictification Theorems I.5.4.6 and I.5.4.7.

Chapter I.9: Matrix Permutative Gray Monoids

This chapter proves variations of Theorem I.8.15.4 when $\mathsf{C}$ satisfies additional strictness conditions. Theorem I.9.3.16 says that, for a *strict* symmetric bimonoidal category $\mathsf{C}$ as in Definition I.9.1.1, $\mathsf{Mat}^{\mathsf{C}}$ is not just a symmetric monoidal bicategory, but is, in fact, a permutative Gray monoid. This means that $\mathsf{Mat}^{\mathsf{C}}$ has an underlying 2-category instead of a bicategory, and it is a monoid with respect to the Gray tensor product for 2-categories. Its braiding is now a symmetry for the Gray monoid structure. Theorem I.9.4.2 says that if, in addition, the multiplicative symmetry $\xi^{\otimes}$ in $\mathsf{C}$ is the identity, then $\mathsf{Mat}^{\mathsf{C}}$ is a permutative 2-category.

**Part II.1: Braided Bimonoidal Categories**

Chapter II.1: Preliminaries on Braided Structures

To prepare for the rest of Part II.1, this chapter discusses the braid groups and braided monoidal categories. First it defines the braid groups and proves some useful properties for sum braids and block braids. Then it reviews braided monoidal categories and proves some basic properties, including two manifestations of the third Reidemeister move. Next it proves in detail that the Drinfeld center of a monoidal category is a braided monoidal category and that the symmetric center of a braided monoidal category is a symmetric monoidal category. Then it recalls the Joyal-Street Coherence Theorem II.1.6.3 for braided monoidal categories.

Chapter II.2: Braided Bimonoidal Categories

This chapter defines braided bimonoidal categories. They are defined using 12 of the 24 Laplaza axioms of a symmetric bimonoidal category, together with two additional axioms that are variants of the only two Laplaza axioms involving the braiding $\xi^{\otimes}$. In a symmetric bimonoidal category, each of these two variant axioms is equivalent to the original Laplaza axiom. This is reminiscent of the fact that a braided monoidal category has two hexagon axioms, which are equivalent to each other in a symmetric monoidal category. A *tight* braided bimonoidal category—that is, one with invertible distributivity morphisms $\delta^l$ and $\delta^r$—is equivalent to a BD category in the sense of Blass and Gurevich [**BG20a**]. The first main observation in this chapter is Theorem II.2.2.1, which says that each braided bimonoidal category satisfies all 24 Laplaza axioms. Therefore, a symmetric bimonoidal category is precisely a braided bimonoidal category whose braiding satisfies the symmetry axiom. The second main result in this chapter says that an abelian category with a compatible (symmetric/braided) monoidal structure is a tight (symmetric/braided) bimonoidal category. The additive structure comes from the abelian structure, and the multiplicative structure comes from the monoidal structure. The braided case of this result is due to Blass and Gurevich [**BG20a**].

### Chapter II.3: Applications to Quantum Groups and Topological Quantum Computation

This chapter shows that braided bimonoidal categories arise naturally in quantum groups and topological quantum computation (TQC). The first main observation is Theorem II.3.2.19. It says that for a (symmetric/braided) bialgebra $A$, the category $\mathsf{Mod}(A)$ of left $A$-modules, equipped with the usual direct sum and tensor product, is a tight (symmetric/braided) bimonoidal category. This is an extension of the important fact in quantum group theory that, for a braided bialgebra $A$, $\mathsf{Mod}(A)$ is a braided monoidal category. Next we prove in detail that Fibonacci anyons and Ising anyons, which are two central models in TQC, are both tight braided bimonoidal categories. In each case, the additive structure comes from an abelian category structure, and the multiplicative structure comes from the fusion rules of anyons.

### Chapter II.4: Bimonoidal Centers

This chapter generalizes the Drinfeld center of a monoidal category and the symmetric center of a braided monoidal category (Sections II.1.4 and II.1.5) to the bimonoidal setting. Generalizing the Drinfeld center, Theorem II.4.4.3 says that, for each tight bimonoidal category $\mathsf{C}$, the bimonoidal Drinfeld center $\overline{\mathsf{C}}^{\mathsf{bi}}$ is a tight braided bimonoidal category. Tightness is required for this theorem because the invertibility of $\delta^l$ and $\delta^r$ is used in the construction of $\overline{\mathsf{C}}^{\mathsf{bi}}$. The proof of this theorem is another good illustration of the axioms of a braided bimonoidal category, since we will use all 24 Laplaza axioms and the two variant axioms in the braided case. Generalizing the symmetric center, Theorem II.4.5.3 says that, for each braided bimonoidal category $\mathsf{C}$, the bimonoidal symmetric center $\mathsf{C}^{\mathsf{sym}}$ is a symmetric bimonoidal category.

### Chapter II.5: Coherence of Braided Bimonoidal Categories

This chapter proves the Coherence Theorem II.5.4.4 for braided bimonoidal categories that satisfy a monomorphism assumption. As in the symmetric case (Theorems I.3.9.1 and I.4.4.3), the monomorphism assumption in Theorem II.5.4.4 is automatically satisfied if tightness is assumed. This theorem is the braided analogue of Laplaza's Second Coherence Theorem I.4.4.3. It uses a braided version $\mathcal{D}^{\mathsf{br}}$ of the distortion category that involves the symmetric groups and the braid groups to keep track of, respectively, the additive symmetry $\xi^{\oplus}$ and the braiding $\xi^{\otimes}$. Reminiscent of the Joyal-Street Coherence Theorem II.1.6.3 for braided monoidal categories, Theorem II.5.4.4 says that, if two paths have the same image in the braided distortion category $\mathcal{D}^{\mathsf{br}}$, then they have the same value in the braided bimonoidal category in question. This condition of having the same image in $\mathcal{D}^{\mathsf{br}}$ is very much checkable in practice. In fact, the proofs of the main results in Chapters II.6, II.7, and II.8 all use Theorem II.5.4.4 and involve checking this condition many times. In [**BG20a**], Blass and Gurevich conjectured the existence of a coherence theorem for their BD categories, which are equivalent to our tight braided bimonoidal categories. Theorem II.5.4.4 confirms the Blass-Gurevich Conjecture in the form of commutative formal diagrams.

Chapter II.6: Strictification of Tight Braided Bimonoidal Categories

This chapter proves two Strictification Theorems II.6.3.6 and II.6.3.7 for tight braided bimonoidal categories. As in the symmetric case (Theorems I.5.4.6 and I.5.4.7), strictification requires the tightness assumption because the construction of the strictification uses the invertibility of the distributivity morphisms $\delta^l$ and $\delta^r$. A *right permbraided category* is a tight braided bimonoidal category with both the additive and the multiplicative structures strict monoidal, and with identities for the structure morphisms $\lambda^{\bullet}$, $\rho^{\bullet}$, $\delta^r$, $\zeta^{\otimes}_{-,0}$, and $\zeta^{\otimes}_{0,-}$. Theorem II.6.3.6 says that each tight braided bimonoidal category is adjoint equivalent to a right permbraided category via strong braided bimonoidal functors. Theorem II.6.3.7 is the analogue that strictifies each tight braided bimonoidal category to a *left* permbraided category, in which $\delta^l$, instead of $\delta^r$, is the identity. Theorems II.6.3.6 and II.6.3.7 are two further positive answers to the Blass-Gurevich Conjecture [**BG20a**] in the form of strictification.

Chapter II.7: The Braided Baez Conjecture

This chapter proves the braided version of Baez's Conjecture. Section II.7.1 defines the 2-category $\mathsf{Bi}_{\mathsf{r}}^{\mathsf{fbr}}$ with *flat* small braided bimonoidal categories as objects. As in the symmetric case, flatness (Definition II.5.4.5) is much weaker than tightness, and it guarantees that the Braided Bimonoidal Coherence Theorem II.5.4.4 is applicable. The first version of the Braided Baez Conjecture (Theorem II.7.3.4) says that the finite ordinal category $\Sigma$ is a lax bicolimit of the 2-functor $\varnothing \longrightarrow \mathsf{Bi}_{\mathsf{r}}^{\mathsf{fbr}}$. Another version is Theorem II.7.3.6, which says that the variant $\Sigma'$ of $\Sigma$ is also such a lax bicolimit. Also like the symmetric case, the proofs of the Braided Baez Conjecture do *not* use the Strictification Theorems II.6.3.6 and II.6.3.7. This allows us to use flat small braided bimonoidal categories in the 2-category $\mathsf{Bi}_{\mathsf{r}}^{\mathsf{fbr}}$, instead of the smaller class of tight ones. The reader may wonder why the finite ordinal category $\Sigma$ and its variant $\Sigma'$ are bi-initial in both the symmetric case (Theorems I.7.8.1 and I.7.8.3) and the braided case. This is analogous to the fact that the ring of integers is initial in both the category of rings and the category of commutative rings.

Chapter II.8: Monoidal Bicategorification

The main Theorem II.8.4.7 in this chapter says that, for each tight braided bimonoidal category $\mathsf{C}$, the matrix construction $\mathsf{Mat}^{\mathsf{C}}$ is a monoidal bicategory. While most of the definitions for $\mathsf{Mat}^{\mathsf{C}}$ are the same as in the symmetric case in Chapter I.8, there are two subtleties. First, in the current braided case, the lax functoriality constraint $\boxtimes^2$ of the monoidal composition $\boxtimes$ in $\mathsf{Mat}^{\mathsf{C}}$ has two additional conditions about the braided distortions of the two paths involved; see (II.8.2.15) and (II.8.2.16). These conditions about the braided distortions are necessary because a braid is not determined by its underlying permutation, and the braided distortion category $\mathcal{D}^{\mathsf{br}}$ involves the braid groups. The second subtle point is that, even if $\mathsf{C}$ is a tight braided bimonoidal category, the monoidal bicategory $\mathsf{Mat}^{\mathsf{C}}$ does *not* seem to have any reasonable braided monoidal bicategory structure in general. We will explain this point in more detail near the end of Section II.8.4. The difficulty once again comes from the fact that the braided distortion category $\mathcal{D}^{\mathsf{br}}$ involves the braid groups, and a braid with an identity underlying permutation is usually not the identity braid.

**Part II.2: $E_n$-Monoidal Categories**

Chapter II.9: Ring, Bipermutative, and Braided Ring Categories

This chapter discusses ring and bipermutative categories in the sense of Elmendorf-Mandell and the braided version. The main difference between these categorical notions and their bimonoidal counterparts in Parts I.1 and II.1 is that ring categories have generally non-invertible *factorization morphisms*

$$(A \otimes C) \oplus (B \otimes C) \xrightarrow{\partial^l_{A,B,C}} (A \oplus B) \otimes C$$
$$(A \otimes B) \oplus (A \otimes C) \xrightarrow{\partial^r_{A,B,C}} A \otimes (B \oplus C)$$

that go in the opposite direction as the distributivity morphisms $\delta^r$ and $\delta^l$. Ring categories with invertible factorization morphisms are special cases of tight bimonoidal categories, so the latter's strictification theorems in Chapter I.5 also apply to such ring categories. The bipermutative and braided analogues are also true. Similar to the endomorphism rig of a commutative monoid, each small permutative category $\mathsf{C}$ yields an endomorphism ring category $\mathsf{Perm}^{\mathsf{su}}(\mathsf{C};\mathsf{C})$. Similar to the reduction of Laplaza's axioms in symmetric bimonoidal categories in Section I.2.2 and the braided version in Theorem II.2.2.1, about half of the ring category axioms are redundant in a bipermutative category and a braided ring category. This is an extension of an observation in [**EM06**, Fig. 1]. Moreover, the Drinfeld center and the symmetric center have natural analogues for these ring-like categories. As we will discuss in Chapters III.11 and III.12, the Elmendorf-Mandell $K$-theory multifunctor sends small ring, braided ring, and bipermutative categories to, respectively, strict ring, $E_2$-, and $E_\infty$-symmetric spectra. The strict ring and $E_\infty$ cases are due to Elmendorf-Mandell [**EM06, EM09**], and the $E_2$ case is new.

Chapter II.10: Iterated and $E_n$-Monoidal Categories

Keeping in mind that the ring-like categories in Chapter II.9 correspond to $E_n$-symmetric spectra for $n \in \{1, 2, \infty\}$ via algebriac $K$-theory, this chapter discusses the categorical structure for the general $E_n$ cases. An $n$-fold monoidal category in the sense of [**BFSV03**] has $n$ monoidal structures $\otimes_1, \ldots, \otimes_n$ that are strictly associative and unital and interact via the exchange natural transformations

$$(A \otimes_j B) \otimes_i (C \otimes_j D) \xrightarrow{\eta^{i,j}_{A,B,C,D}} (A \otimes_i C) \otimes_j (B \otimes_i D)$$

for $1 \leq i < j \leq n$. Monoids in the monoidal category of small $n$-fold monoidal categories are precisely small $(n+1)$-fold monoidal categories. We introduce the notion of an *$E_n$-monoidal category* as a permutative category $(\mathsf{C}, \oplus)$ equipped with an $n$-fold monoidal structure $\{\otimes_i, \eta^{i,j}\}$ and factorization morphisms $\{\partial^{l,i}, \partial^{r,i}\}$ for each monoidal structure $\otimes_i$, such that (i) each $(\oplus, \otimes_i, \partial^{l,i}, \partial^{r,i})$ is a ring category and (ii) several axioms relating $\eta^{i,j}$, $\partial^{l,i}$, and $\partial^{r,i}$ hold. Ring categories are $E_1$-monoidal categories. Braided ring categories and bipermutative categories are special cases of, respectively, $E_2$- and $E_n$-monoidal categories for $n \geq 2$. Moreover, each small category generates a free $E_n$-monoidal category. In Chapter III.13, we will show that the Elmendorf-Mandell $K$-theory of a small $E_n$-monoidal category is an $E_n$-symmetric spectrum for $n \geq 1$.

**Part III.1: Enriched Monoidal Categories and Multicategories**

Chapter III.1: Enriched Monoidal Categories

This chapter gives the basic definitions and properties for enriched monoidal categories, including plain, braided, and symmetric variants. Definition III.1.4.25 describes 2-categories of each, with 1- and 2-cells given by appropriately monoidal enriched functors and natural transformations, respectively. For our applications to $K$-theory in Part III.2, the enriching category $\mathsf{V}$ is symmetric monoidal closed. However, our treatment in this chapter addresses the more general case that $\mathsf{V}$ is merely monoidal, with additional assumptions about braided or symmetric monoidal structure stated as necessary.

Section III.1.5 discusses the important special case $\mathsf{V} = \mathsf{Cat}$, the category of small categories with its Cartesian product. Explanation III.1.5.3 describes how the monoidal $\mathsf{V}$-categories in this case are strict versions of monoidal bicategories. The braided and symmetric cases are similarly compared.

Chapter III.2: Change of Enrichment

This chapter describes change of enriching category induced by a symmetric monoidal functor, showing that monoidal structures are preserved. Sections III.2.1 through III.2.4 give a thorough treatment of 2-functoriality results. As an application, Corollary III.2.4.17 shows that taking underlying categories gives a 2-functor from small monoidal $\mathsf{V}$-categories, $\mathsf{V}$-functors, and $\mathsf{V}$-natural transformations to ordinary monoidal categories, functors, and natural transformations. Similar statements hold for the braided and symmetric cases.

A partial reverse of Corollary III.2.4.17 is given in Theorem III.2.5.1. The theorem shows that, for given $\mathsf{V}$-enriched data, various enriched monoidal axioms are satisfied if and only if the corresponding monoidal axioms for the underlying data are satisfied. This provides a mechanism to lift ordinary monoidal structures to enriched monoidal structures.

Sections III.2.5 and III.2.6 apply Theorem III.2.5.1 to lift coherence and strictification results for ordinary monoidal, braided, and symmetric monoidal categories to their enriched counterparts. The Enriched Monoidal Coherence Theorem III.2.5.6 and Enriched Epstein's Coherence Theorem III.2.5.8 play a significant role in subsequent chapters.

Chapter III.3: Self-Enrichment and Enriched Yoneda

This chapter restricts to the case that $\mathsf{V}$ is a symmetric monoidal closed category. Theorem III.3.3.2 shows, via Theorem III.2.5.1, that the canonical enrichment of $\mathsf{V}$ over itself is symmetric monoidal as a $\mathsf{V}$-category. The next several sections develop the theory of $\mathsf{V}$-enriched co/ends followed by the $\mathsf{V}$-Yoneda Lemma (Theorem III.3.6.9) and an equivalent form called the $\mathsf{V}$-Yoneda Density Theorem III.3.7.8. These are applied to develop the Day convolution and internal hom for enriched diagram categories (Theorem III.3.7.22). The remainder of the chapter discusses additional theory of enriched diagram categories and tensor/cotensor structures that will be important for the development of enriched $K$-theory functors in Part III.2.

Chapter III.4: Pointed Objects, Smash Products, and Pointed Homs

This chapter gives the definitions and properties of smash products and pointed homs. These will be used throughout Part III.2, and the smash product of pointed multicategories, developed in Chapter III.5, will be particularly important.

Section III.4.3 uses the Day convolution and internal hom to develop symmetric monoidal closed structure for pointed diagram categories. The results are summarized in Theorem III.4.3.37. Applications of this material appear in Chapters III.8, III.9, and III.10, where the Segal and Elmendorf-Mandell $K$-theory constructions are given via certain pointed diagram categories.

Chapter III.5: Multicategories

This chapter gives relevant background on multicategories, multifunctors, and multinatural transformations. Theorem III.5.5.14 shows that the category of small multicategories is complete and cocomplete. The Boardman-Vogt tensor product of multicategories, and the associated smash product for pointed multicategories, are developed in Section III.5.6. The corresponding internal hom and its pointed variant are developed in Section III.5.7.

Chapter III.6: Enriched Multicategories

This chapter develops basic definitions and properties for enriched multicategories. One of our important applications, developed in Section III.6.3, is the enriched multicategory associated to an enriched symmetric monoidal category. Our first use of this is in Section III.6.4 where we describe the Cat-enriched multicategory structure on Multicat, the category of small multicategories. It is induced by showing that the tensor product makes Multicat symmetric monoidal as a Cat-enriched category (Theorem III.6.4.3). The pointed variant, with the smash product of small pointed multicategories, is given in Theorem III.6.4.4 and will be essential for Part III.2.

Sections III.6.5 and III.6.6 cover our second important application of enriched multicategories. The category $\mathsf{PermCat}^{\mathsf{su}}$, consisting of small permutative categories and strictly unital symmetric monoidal functors, has a Cat-enriched multicategory structure given by multilinear functors and multilinear transformations (Definitions III.6.5.4 and III.6.5.11). Propositions III.6.5.10 and III.6.5.13 show that this Cat-enriched multicategory structure is induced from that of small pointed multicategories and their smash product. Section III.6.6 gives a second, direct proof of the Cat-enriched multicategory axioms.

**Part III.2: Algebraic *K*-Theory**

Chapter III.7: Homotopy Theory Background

This chapter gives relevant background from homotopy theory. Sections III.7.1 and III.7.2 introduce simplicial sets and simplicial homotopy, along with the nerve and geometric realization functors. The category of symmetric spectra, with its symmetric monoidal closed structure, is presented in Sections III.7.3 through III.7.6. Then, Sections III.7.7 and III.7.8 give a short review of Quillen model categories and a number of key examples.

Chapter III.8: Segal $K$-Theory of Permutative Categories

This chapter presents the $K$-theory functor $\mathsf{K}^{\mathsf{Se}}$ due to Segal [**Seg74**]. Its inputs are small permutative categories and its outputs are symmetric spectra. Section III.8.3 describes the key construction as given by Segal. Sections III.8.4 and III.8.5 describe an equivalent construction that compares more readily with the $K$-theory multifunctor of Elmendorf-Mandell, $\mathsf{K}^{\mathsf{EM}}$.

Chapter III.9: Categories of $\mathcal{G}_*$-Objects

This chapter is the first of two that replace the Segal $K$-theory functor with a simplicially-enriched multifunctor due to Elmendorf-Mandell [**EM06, EM09**]. This chapter focuses on the replacement of $\Gamma$-categories and $\Gamma$-simplicial sets with pointed diagrams out of a larger indexing category $\mathcal{G}$. The construction of symmetric spectra from such diagram categories is given in Section III.9.3 and is denoted $\mathsf{K}^{\mathcal{G}}$. Sections III.9.2 and III.9.4 use the material from Part III.1 to explain that the new diagram categories and the new functor $\mathsf{K}^{\mathcal{G}}$ are symmetric monoidal, in the enriched sense of Chapter III.1, over the category of pointed simplicial sets.

Chapter III.10: Elmendorf-Mandell $K$-Theory of Permutative Categories

This chapter is the second of two that replace the Segal $K$-theory functor with a simplicially enriched multifunctor due to Elmendorf-Mandell [**EM06, EM09**]. This chapter focuses on the construction of $\mathcal{G}_*$-categories from small permutative categories, replacing Segal's construction of $\Gamma$-categories from the same. Additional material from Part III.1 is used throughout the chapter to explain that the multi/categories and multi/functors are enriched either in the symmetric monoidal sense of Chapter III.1 or in the multicategorical and multifunctorial sense of Chapter III.6. Section III.10.6 contains the proof that the Segal and Elmendorf-Mandell $K$-theory symmetric spectra associated to a small permutative category $\mathsf{C}$ are level equivalent (Theorem III.10.6.10). Because $\mathsf{K}^{\mathsf{EM}}$ is an enriched multifunctor, it preserves operad actions. We state this result as Theorem III.10.3.33 and apply it in Chapters III.11, III.12, and III.13.

Chapter III.11: $K$-Theory of Ring and Bipermutative Categories

This is the first of three chapters that contain algebraic $K$-theory applications of the ring-like categories in Part II.2. The main $K$-theory results in this chapter, Corollaries III.11.3.16 and III.11.6.12, are from [**EM06, EM09**], and they are the $E_1$ and the $E_\infty$ cases. These results state that the Elmendorf-Mandell $K$-theory multifunctor $\mathsf{K}^{\mathsf{EM}}$ sends (i) small ring categories to strict ring symmetric spectra and (ii) small bipermutative categories to $E_\infty$-symmetric spectra. They are obtained by combining the multifunctor $\mathsf{K}^{\mathsf{EM}}$ and the fact that the associative operad and the Barratt-Eccles operad parametrize, respectively, ring and bipermutative category structures on small permutative categories. Since the associative operad has monoids as algebras and the Barratt-Eccles operad is an $E_\infty$-operad, the $K$-theory results follow.

Chapter III.12: $K$-Theory of Braided Ring Categories

This chapter contains the $E_2$ analogues of the results in Chapter III.11. The first part of this chapter discusses the braid operad $\mathsf{Br}$, which generalizes the Barratt-Eccles operad. This is a categorical $E_2$-operad (Theorem III.12.2.4) whose algebras in $\mathsf{Cat}$ are small braided strict monoidal categories (Proposition III.12.3.22). The main categorical input is Theorem III.12.4.5, which says that $\mathsf{Br}$ parametrizes braided ring category structures, as in Chapter II.9, on small permutative categories. Applying the $K$-theory multifunctor $\mathsf{K}^{\mathsf{EM}}$, it follows that $\mathsf{K}^{\mathsf{EM}}$ sends small braided ring categories to $E_2$-symmetric spectra (Corollary III.12.5.3). The $K$-theory result, Corollary III.12.5.3, and the main categorical input, Theorem III.12.4.5, are new results.

Chapter III.13: $K$-Theory of $E_n$-Monoidal Categories

This chapter contains the general $E_n$ analogues for $n \geq 1$ of the categorical and $K$-theory results in Chapters III.11 and III.12. The first part of this chapter discusses the $n$-fold monoidal category operad $\mathsf{Mon}^n$. This is a categorical $E_n$-operad (Theorem III.13.2.1) whose algebras in $\mathsf{Cat}$ are small $n$-fold monoidal categories (Proposition III.13.3.18) as in Chapter II.10. The main categorical input is Theorem III.13.4.12, which says that $\mathsf{Mon}^n$ parametrizes $E_n$-monoidal category structures on small permutative categories. Applying the $K$-theory multifunctor $\mathsf{K}^{\mathsf{EM}}$, it follows that $\mathsf{K}^{\mathsf{EM}}$ sends small $E_n$-monoidal categories to $E_n$-symmetric spectra for $n \geq 1$ (Corollary III.13.5.2). As in Chapter III.12, the $K$-theory result, Corollary III.13.5.2, and the main categorical input, Theorem III.13.4.12, are new results.

Appendix III.A: Open Questions

This chapter discusses open questions related to the topics of this work. We encourage the reader to read these open questions at any time and use them as additional motivation for the main text.

# Bibliography and Indices

APPENDIX A

# Open Questions

*"I enjoy questions that seem honest, even when they admit or reveal confusion, in preference to questions that appear designed to project sophistication."*

– Bill Thurston, MathOverflow user profile

In this chapter, we discuss open questions related to the topics of this work. These open questions provide additional motivation for the main text.

## A.1. Bimonoidal Categories

The following questions are about bimonoidal, symmetric bimonoidal, and braided bimonoidal categories in Definitions I.2.1.2 and II.2.1.29.

**Question A.1.1** (Functoriality of the Matrix Construction)**.** In Theorem I.8.15.4, we showed that, for each tight symmetric bimonoidal category $\mathsf{C}$, the matrix construction $\mathsf{Mat}^{\mathsf{C}}$ is a symmetric monoidal bicategory. Denote by $\mathsf{Bi}^{\mathsf{tsy}}$ the full sub-2-category of the 2-category $\mathsf{Bi}^{\mathsf{sy}}$ in Proposition I.7.1.7, with small *tight* symmetric bimonoidal categories as objects. Regard $\mathsf{Bi}^{\mathsf{tsy}}$ as a tricategory with only identity 3-cells. It is claimed in [**SP∞**] that small symmetric monoidal bicategories are the objects of a tricategory, denoted by $\mathsf{SMB}$.

- Extend the assignment

$$\mathsf{C} \longmapsto \mathsf{Mat}^{\mathsf{C}}$$

  to a trifunctor

$$\mathsf{Bi}^{\mathsf{tsy}} \longrightarrow \mathsf{SMB}.$$

Among other things, one should carefully verify the tricategory axioms for $\mathsf{SMB}$. For a discussion of tricategories and a detailed verification of the tricategory of small bicategories, the reader is referred to [**JY21**, Ch. 11]. ◇

**Question A.1.2** (Bimonoidal Bicategories)**.** Taking the categorification from (commutative) rigs to (symmetric) bimonoidal categories one step further, we could ask about two different monoidal structures, $\boxplus$ and $\boxtimes$, on a bicategory, with $\boxplus$ symmetric.

- Define such a *(braided/sylleptic/symmetric) bimonoidal bicategory*, generalizing the (braided/sylleptic/symmetric) monoidal bicategories in Sections I.6.4 and I.6.5.
- For a tight symmetric bimonoidal category $\mathsf{C}$, prove that the symmetric monoidal bicategory $\mathsf{Mat}^{\mathsf{C}}$ extends to a symmetric bimonoidal bicategory.
- For a tight braided bimonoidal category $\mathsf{C}$, prove that the monoidal bicategory $\mathsf{Mat}^{\mathsf{C}}$ extends to a bimonoidal bicategory.

More precisely, in Theorem I.8.15.4, the monoidal composition $\boxtimes$ in the symmetric monoidal bicategory $\mathsf{Mat}^{\mathsf{C}}$ involves the matrix tensor product in (I.8.6.3). There should be another symmetric monoidal bicategory structure on the matrix bicategory $\mathsf{Mat}^{\mathsf{C}}$ in Theorem I.8.4.12, whose monoidal composition uses the matrix direct sum as in Example I.2.5.9. These two symmetric monoidal structures on $\mathsf{Mat}^{\mathsf{C}}$ should make it into a symmetric bimonoidal bicategory. A similar discussion applies in the braided case, with Theorem II.8.4.7 showing that $\mathsf{Mat}^{\mathsf{C}}$ is a monoidal bicategory.

- Extend the tricategory SMB (Question A.1.1) to a tricategory SBB with small symmetric bimonoidal bicategories as objects.
- Construct a tricategory BB with small bimonoidal bicategories as objects.
- Extend the assignment $\mathsf{C} \longmapsto \mathsf{Mat}^{\mathsf{C}}$ in
  - Theorem I.8.15.4 in the symmetric case to a trifunctor

$$\mathsf{Bi}^{\mathsf{tsy}} \longrightarrow \mathsf{SBB}$$

    and
  - Theorem II.8.4.7 in the braided case to an analogous trifunctor with codomain BB. ◇

**Question A.1.3** (Bimonoidal Bicategorical Centers)**.** Extend the bimonoidal centers in Theorems II.4.4.3 and II.4.5.3 to the bimonoidal bicategorical setting in Question A.1.2. In other words:

- Extend the bimonoidal Drinfeld center in Theorem II.4.4.3 to a bimonoidal bicategory and show that it is a braided bimonoidal bicategory.
- Show that the center of a braided bimonoidal bicategory is a sylleptic bimonoidal bicategory.
- Extend the bimonoidal symmetric center in Theorem II.4.5.3 to a sylleptic bimonoidal bicategory and show that it is a symmetric bimonoidal bicategory.

There are simpler centers of Gray monoids **[BN96]**, braided monoidal 2-categories, and sylleptic monoidal 2-categories **[Cra98]**. As in Theorem II.4.4.3, a tightness assumption is likely necessary for some of these center constructions. ◇

**Question A.1.4** (Gray Rings and Bipermutative Gray Monoids)**.** Recall from Sections I.6.6 and I.6.7 that a Gray monoid is a 2-category equipped with a monoid structure under the Gray tensor product. A permutative Gray monoid is a Gray monoid equipped with a compatible Gray symmetry. Symmetric monoidal bicategories can be strictified to permutative Gray monoids; see **[GJO17b]** and the discussion near the end of Section I.6.7.

- Define a *Gray ring* and a *bipermutative Gray monoid* that are analogous to, respectively, a right rigid bimonoidal category (Definition I.5.5.8) and a right bipermutative category (Definition I.2.5.2).
- Along the lines of Theorem I.5.5.11, prove a strictification result from bimonoidal bicategories (Question A.1.2) to Gray rings.
- Along the lines of Theorem I.5.4.6, prove a strictification result from symmetric bimonoidal bicategories to bipermutative Gray monoids.

A bipermutative Gray monoid should be a 2-category equipped with two compatible permutative Gray monoid structures, $\boxplus$ and $\boxtimes$, that interact via distributivity.

In a Gray ring, $\boxtimes$ is a Gray monoid structure that is not assumed to be permutative. The following table summaries these (conjectural) concepts.

| lax structure | strict structure | strictification |
|---|---|---|
| bicategories | 2-categories | [**JY21**, 8.4.1] |
| monoidal bicategories | Gray monoids | [**GPS95, Gur13**] |
| symmetric monoidal bicategories | permutative Gray monoids | [**GJO17b**] |
| bimonoidal bicategories | Gray rings | conjecture |
| symmetric bimonoidal bicategories | bipermutative Gray monoids | conjecture |

In each row, the left column can be strictified to the middle column. ◇

**Question A.1.5** (Horizontal Bicategories of Double Categories)**.** A number of bicategories, including those of spans and bimodules, are the horizontal bicategories of some double categories.

- For a tight bimonoidal category $\mathsf{C}$, is the matrix bicategory $\mathsf{Mat}^{\mathsf{C}}$ in Theorem I.8.4.12 the horizontal bicategory of a double category?
- If so, does the symmetric monoidal bicategory in Theorem I.8.15.4 arise from a symmetric monoidal structure on the double category?

See [**HS∞**] and [**JY21**, 12.3 and 12.4] for a discussion of (monoidal) double categories and their horizontal bicategories. ◇

**Question A.1.6** (Strict Symmetric Bimonoidal Categories)**.** The Strictification Theorems I.5.4.6 and I.5.4.7 state that each tight symmetric bimonoidal category is adjoint equivalent to a right bipermutative category and a left bipermutative category.

- Is there an analogue that strictifies a tight symmetric bimonoidal category to a *strict* symmetric bimonoidal category as in Definition I.9.1.1?

By Theorems I.5.4.6 and I.5.4.7, one may start with a right or left bipermutative category. ◇

**Question A.1.7** (Braided Sheet Diagrams)**.** String diagrams are graphical reasoning tools in monoidal categories [**JS91a, Sel11**]. Sheet diagrams [**CDH∞**], which we mentioned in Notes I.2.7.5 and I.7.9.2 and Example I.3.10.9, are their analogues for tight bimonoidal categories.

- Develop sheet diagrams for
  - symmetric bimonoidal categories (Definition I.2.1.2) and
  - braided bimonoidal categories (Definition II.2.1.29).

This is, in fact, a coherence question with several parts. More precisely, a *bimonoidal signature* $S$ consists of (i) a set of generating objects and (ii) a set of generating morphisms, each with (co)domain in the free $\{\oplus,\otimes\}$-algebra $S^{\mathsf{fr}}$ (Definition I.3.1.2). Given a bimonoidal signature $S$, one first defines the appropriate braided bimonoidal sheet diagrams and topological deformations corresponding to the axioms of a braided bimonoidal category. Then one constructs a braided bimonoidal category $S'$ with object set $S^{\mathsf{fr}}$ and, as morphisms, braided bimonoidal sheet diagrams modulo topological deformations. Finally, one proves that $S'$ is braided bimonoidally equivalent to the free braided bimonoidal category on $S$. A similar discussion applies in the symmetric case.

In view of the results in Chapter II.3, braided bimonoidal sheet diagrams can be used as graphical reasoning tools in (i) quantum group theory and (ii) the Fibonacci and Ising anyons in topological quantum computation. Sheet diagrams for tight bimonoidal categories in [**CDH∞**] involve the symmetric monoidal string diagrams in [**JS91a**]. Braided bimonoidal sheet diagrams will likely involve both the symmetric monoidal (for the additive structure $\oplus$) and the braided monoidal (for the multiplicative structure $\otimes$) string diagrams in [**JS91a**, Ch. 2–3]. The Braided Bimonoidal Coherence Theorem II.5.4.4 will be needed to check the axioms in Definition II.2.1.29 for braided bimonoidal sheet diagrams. Similarly, Laplaza's Coherence Theorems I.3.9.1 and I.4.4.3 will be needed to check the axioms in Definition I.2.1.2 for symmetric bimonoidal sheet diagrams.

The distributivity morphisms $\delta^l$ and $\delta^r$ in Definitions I.2.1.2 and II.2.1.29 are not invertible in general.

- Is it possible to replace the tightness assumption—that is, the invertibility of $\delta^l$ and $\delta^r$—in the sheet diagrams in [**CDH∞**] with flatness in Definitions I.3.9.9 and II.5.4.5?

Related to Section I.7.9, to replace tightness with the much weaker assumption of flatness, one would need to work directly with a flat (symmetric/braided) bimonoidal category and *avoid* using the Strictification Theorems I.5.4.6, I.5.4.7, I.5.5.11, I.5.5.12, II.6.3.6, and II.6.3.7. The reason is that each of these theorems requires the tightness assumption. See Question A.2.8 for further problems about sheet diagrams. ⋄

## A.2. $E_n$-Monoidal Categories

The following questions are about the $E_n$-monoidal categories in Part II.2.

**Question A.2.1** (Coherence of $E_n$-Monoidal Categories)**.**

- Prove coherence theorems for ring categories (Definition II.9.1.2) along the lines of Theorems I.3.10.7 and I.4.5.8. Each such coherence theorem should say that any reasonable formal diagram in a ring category involving
$$(\oplus, 0, \xi^{\oplus}, \otimes, \mathbb{1}, \partial^l, \partial^r)$$
is commutative, with an assumption on either the common domain or the two paths.
- Prove coherence theorems for bipermutative categories (Definition II.9.3.2) along the lines of Theorems I.3.9.1 and I.4.4.3.
- Prove a coherence theorem for braided ring categories (Definition II.9.5.1) along the lines of Theorem II.5.4.4.
- More generally, prove a coherence theorem for $E_n$-monoidal categories (Definition II.10.7.2) along the lines of Theorem II.5.4.4. The coherence theorems for $n$-fold monoidal categories (Theorem II.10.6.8) should be relevant.

As in Theorems I.3.9.1 and I.4.4.3, one may need to assume a monomorphism or an epimorphism condition on the factorization morphisms $\partial^l$ and $\partial^r$ ⋄

**Question A.2.2** ($n$-Monoidal Categories)**.** In an $n$-fold monoidal category (Definition II.10.1.1), each monoidal structure $\otimes_i$ is strictly associative with a common strict unit $\mathbb{1}$. There is a more general concept called an *$n$-monoidal category* in [**AM10**, Def. 6.1, 7.1, and 7.24]. It allows each monoidal structure $\otimes_i$ to be nonstrict and distinct monoidal units.

- Describe the free $n$-monoidal category of a small category, along the lines of Proposition II.10.5.9 and Theorem II.10.5.18.
- Generalize the Coherence Theorem II.10.6.8 to $n$-monoidal categories.

In [**AM10**, Section 6.2], it is stated without detail that, in a 2-monoidal category, each formal diagram is commutative. It is stated there that this coherence result can be deduced from the work in [**Lew72**]. So one possible first step in answering these questions would be to prove in detail this coherence result for 2-monoidal categories. ◇

Questions A.2.3 through A.2.7 below are all related to each other.

**Question A.2.3** (Lax $n$-Fold Monoidal Categories)**.** Between an $n$-fold monoidal category and an $n$-monoidal category (Question A.2.2) is a *lax n-fold monoidal category*. The latter allows each monoidal structure $\otimes_i$ to be nonstrict, and it assumes a common nonstrict monoidal unit $\mathbb{1}$. Analogous to Propositions II.10.1.14 and II.10.1.21, general nonstrict braided monoidal categories should be examples of lax 2-fold monoidal categories, and general nonstrict symmetric monoidal categories should be examples of lax $n$-fold monoidal categories for $n \geq 2$.

- Describe the free lax $n$-fold monoidal category of a small category, along the lines of Proposition II.10.5.9 and Theorem II.10.5.18.
- Extend the Coherence Theorem II.10.6.8 to lax $n$-fold monoidal categories.
- Can $n$-monoidal categories be strictified to lax $n$-fold monoidal categories?
- Can lax $n$-fold monoidal categories be strictified to $n$-fold monoidal categories?

There are two other possible variants of $n$-fold monoidal categories. The variant in [**For04**] corresponds to a lax $n$-fold monoidal category with a common *strict* monoidal unit. The variant in [**FSS07**] corresponds to an $n$-monoidal category with generally distinct but strict monoidal units. The following table summaries the strictness assumptions of $n$-fold monoidal categories and its four variants.

| monoidal units | strict $\{\otimes_i\}_{1\leq i\leq n}$ | nonstrict $\{\otimes_i\}_{1\leq i\leq n}$ |
|---|---|---|
| common strict | $n$-fold monoidal | [**For04**] |
| common nonstrict | | lax $n$-fold monoidal |
| distinct strict | | [**FSS07**] |
| distinct nonstrict | | $n$-monoidal [**AM10**] |

The questions below refer to lax $n$-fold monoidal categories as defined above, with $\otimes_i$ generally nonstrict and a common nonstrict monoidal unit. ◇

**Question A.2.4** (Unstable Periodic Table of Weak $n$-Categories)**.** The *periodic table* in [**BD98**] of $k$-tuply monoidal $n$-categories is a guiding principle for defining some versions of weak $n$-categories. In the $n = 1$ column in the periodic table, the values $k = 0, 1, 2$, and $\geq 3$ correspond to, respectively, categories, monoidal categories, braided monoidal categories, and symmetric monoidal categories. On the other hand, by Propositions II.10.1.14 and II.10.1.21, braided strict monoidal categories and permutative categories are special examples of, respectively, 2-fold monoidal categories and $k$-fold monoidal categories for $k \geq 2$. Proposition II.10.2.8 has examples of 2-fold monoidal categories that are not braided strict monoidal categories. The $k$-fold monoidal category operad $\mathsf{Mon}^k$ (Proposition 13.1.20) is an $E_k$-operad that parametrizes $k$-fold monoidal categories (Theorem 13.2.1 and Proposition 13.3.18).

- Construct an *unstable* periodic table in which the $n = 1$ column consists of lax $k$-fold monoidal categories (Question A.2.3) for $k \geq 1$.

The $n = 1$ and $n = 2$ columns of the unstable periodic table should look like this:

| | $n = 1$ | $n = 2$ |
|---|---|---|
| $k = 0$ | categories | bicategories |
| $k = 1$ | monoidal categories | monoidal bicategories |
| $k \geq 2$ | lax $k$-fold monoidal categories | $k$-fold monoidal bicategories |

The unstable periodic table does *not* stabilize like the periodic table in **[BD98]**, where the $(n = 1, k \geq 3)$ entries are all symmetric monoidal categories. Moreover, the $(n = 2, k \geq 2)$ column in the unstable periodic table contains yet-to-be-defined *k-fold monoidal bicategories*.

- Prove that braided, sylleptic, and symmetric monoidal bicategories in Section I.6.5 are examples of $k$-fold monoidal bicategories for, respectively, $k = 2, 3$, and $\geq 4$. ◇

This question may be regarded as both (i) a litmus test for the correct definition of a $k$-fold monoidal bicategory and (ii) a conceptual unification of braided, sylleptic, and symmetric monoidal bicategories. Further examples of $k$-fold monoidal bicategories should arise from the matrix construction in Question A.2.7.

**Question A.2.5** (Iterated Gray Monoids)**.** This is a variation of Question A.1.4.

- Analogous to $k$-fold monoidal categories (Definition II.10.1.1) for $k \geq 1$, define the concept of a *k-fold Gray monoid* that satisfies the following statements:
  - A Gray monoid is precisely a 1-fold Gray monoid.
  - A braided monoidal 2-category **[Cra98]** is an example of a 2-fold Gray monoid, analogous to Proposition II.10.1.14.
  - A sylleptic monoidal 2-category **[Cra98]** is an example of a 3-fold Gray monoid.
  - A permutative Gray monoid is an example of a $k$-fold Gray monoid for $k \geq 4$, analogous to Proposition II.10.1.21.
- Prove a strictification theorem from $k$-fold monoidal bicategories (Question A.2.4) to $k$-fold Gray monoids. This should fit into the following table of strictification theorems.

| | | |
|---|---|---|
| bicategories | 2-categories | **[JY21**, 8.4.1**]** |
| monoidal bicategories | Gray monoids | **[GPS95, Gur13]** |
| braided monoidal bicategories | braided monoidal 2-categories | **[Gur11]** |
| $k$-fold monoidal bicategories ($k \geq 1$) | $k$-fold Gray monoids | conjecture |
| symmetric monoidal bicategories | permutative Gray monoids | **[GJO17b]** |

In each row, the left column can be strictified to the middle column. For the symmetric case, see the discussion near the end of Section I.6.7. ◇

**Question A.2.6** (Laplaza $E_n$-Monoidal Categories)**.** The factorization morphisms $\{\partial^{l,i}, \partial^{r,i}\}_{1 \le i \le n}$ in an $E_n$-monoidal category (Definition II.10.7.2) go in the opposite direction as the distributivity morphisms in a bimonoidal category (Definition I.2.1.2). Moreover, in an $E_n$-monoidal category, the monoidal structures $\oplus$ and $\{\otimes_i\}_{i=1}^n$ are all strict.

- Define a *Laplaza $E_n$-monoidal category* with
  - a generally nonstrict additive structure $(\oplus, 0)$,
  - a lax $n$-fold monoidal structure $\big(\{\otimes_i\}_{i=1}^n, \{\eta^{i,j}\}\big)$ (Question A.2.3),
  - for each $1 \le i \le n$, a bimonoidal structure $(\oplus, \otimes_i, \lambda_i^\bullet, \rho_i^\bullet, \delta_i^l, \delta_i^r)$ (Definition I.2.1.2), and
  - appropriate axioms relating the lax $n$-fold monoidal structure and the $n$ bimonoidal structures, analogous to those in Definition II.10.7.2.

Laplaza $E_n$-monoidal categories should contain the following examples:

- An $E_n$-monoidal category with invertible factorization morphisms should be an example of a Laplaza $E_n$-monoidal category, analogous to Theorems II.9.1.15, II.9.3.7, and II.9.5.6.
- A Laplaza $E_1$-monoidal category should be precisely a bimonoidal category (Definition I.2.1.2), analogous to Example II.10.7.13.
- A braided bimonoidal category (Definition II.2.1.29) should be an example of a Laplaza $E_2$-monoidal category with $\otimes_1 = \otimes_2$, analogous to Theorem II.10.8.1.
- A symmetric bimonoidal category should be an example of a Laplaza $E_n$-monoidal category for $n \ge 2$ with $\otimes_1 = \cdots = \otimes_n$, analogous to Theorem II.10.9.1.
- Similar to Theorem II.2.4.22, an abelian category with a compatible lax $n$-fold monoidal structure should be a Laplaza $E_n$-monoidal category.

Moreover:

- Prove a coherence theorem for Laplaza $E_n$-monoidal categories, along the lines of Theorem II.5.4.4. This will certainly involve the coherence theorem for lax $n$-fold monoidal categories in Question A.2.3.
- Prove a strictification theorem for *tight* Laplaza $E_n$-monoidal categories, along the lines of Theorems II.6.3.6 and II.6.3.7. Here *tight* means that all the distributivity morphisms, $\delta_i^l$ and $\delta_i^r$, are natural isomorphisms.
- Is there an analogue of Baez's Conjecture (Theorems I.7.8.1, I.7.8.3, II.7.3.4, and II.7.3.6) for Laplaza $E_n$-monoidal categories? ⋄

**Question A.2.7** (Matrix Construction)**.** In Theorem I.8.4.12 we showed that the matrix construction $\mathsf{Mat}^{\mathsf{C}}$ is a bicategory for each tight bimonoidal category $\mathsf{C}$. Moreover, $\mathsf{Mat}^{\mathsf{C}}$ is (i) a monoidal bicategory if $\mathsf{C}$ is a tight braided bimonoidal category (Theorem II.8.4.7) and (ii) a symmetric monoidal bicategory if $\mathsf{C}$ is a tight symmetric bimonoidal category (Theorem I.8.15.4).

- Show that the matrix construction $\mathsf{Mat}^{\mathsf{C}}$ of a tight Laplaza $E_{k+1}$-monoidal category (Question A.2.6) is a $k$-fold monoidal bicategory.

This question asks for a refinement of the table in the introduction of Chapter II.8 as follows.

| tight – category C | – bicategory $\mathsf{Mat}^{\mathsf{C}}$ | |
|---|---|---|
| bimonoidal | plain | I.8.4.12 |
| braided bimonoidal | monoidal | II.8.4.7 |
| Laplaza $E_{k+1}$-monoidal ($k \geq 1$) | $k$-fold monoidal | conjecture |
| symmetric bimonoidal | symmetric monoidal | I.8.15.4 |

The $k$-fold monoidal bicategories in the conjectural row refer to the $n = 2$ column in the unstable periodic table in Question A.2.4. Proving that $\mathsf{Mat}^{\mathsf{C}}$ is a $k$-fold monoidal bicategory will certainly involve the coherence theorem for Laplaza $E_{k+1}$-monoidal categories in Question A.2.6. The general picture of the table above is that it takes a sum $\oplus$ and a product $\otimes$ to construct the matrix bicategory $\mathsf{Mat}^{\mathsf{C}}$. So any further monoidal structures on the bicategory $\mathsf{Mat}^{\mathsf{C}}$ would have to come from further monoidal structures on C. ◇

**Question A.2.8** (Higher Sheet Diagrams)**.** As in Question A.1.7, develop higher dimensional sheet diagrams for

- $n$-fold monoidal categories (Definition II.10.1.1),
- $n$-monoidal categories (Question A.2.2),
- lax $n$-fold monoidal categories (Question A.2.3),
- $E_n$-monoidal categories (Definition II.10.7.2), and
- Laplaza $E_n$-monoidal categories (Question A.2.6).

As discussed in Question A.1.7, each of these items is a coherence question with several parts. ◇

## A.3. Enriched Monoidal Categories

The following questions are about the concepts in Chapters 1, 2, and 3.

**Question A.3.1** (Enriched Lax $n$-Fold Monoidal Categories)**.** In Section 1.4, with V a braided monoidal category, we defined monoidal, braided monoidal, and symmetric monoidal V-categories, with the latter two assuming that V is symmetric. See Lemma 1.3.23 and Explanation 1.3.25.

- For $n \geq 2$ and V a symmetric monoidal category, extend the lax $n$-fold monoidal categories in Question A.2.3 to the V-enriched setting.
- Extend the results in Chapter 2 to lax $n$-fold monoidal V-categories.

Analogous to Propositions II.10.1.14 and II.10.1.21, braided monoidal V-categories should be examples of lax 2-fold monoidal V-categories, and symmetric monoidal V-categories should be examples of lax $n$-fold monoidal V-categories for $n \geq 2$. Theorem II.10.4.5 says that small $(n+1)$-fold monoidal categories are precisely the monoids in the monoidal category $\mathsf{MCat}^n$ of small $n$-fold monoidal categories and $n$-fold monoidal functors.

- Extend Theorem II.10.4.5 to the V-enriched setting. ◇

**Question A.3.2** (Centers and Enriched Centers)**.** The Drinfeld center of a monoidal category is a braided monoidal category (Theorem II.1.4.27), and the symmetric center of a braided monoidal category is a symmetric monoidal category (Proposition II.1.5.3). Moreover, the Drinfeld center and the symmetric center are generalized to (i) the bimonoidal setting in Theorems II.4.4.3 and II.4.5.3 and (ii) the ring categorical setting in Corollary II.9.6.1 and Theorem II.9.6.4.

- Define a center construction that sends an $n$-fold monoidal category (Definition II.10.1.1) to an $(n+1)$-fold monoidal category.

Theorem II.10.4.5 should be relevant.

- Repeat the previous question for
  - $n$-monoidal categories (Question A.2.2),
  - lax $n$-fold monoidal categories (Question A.2.3),
  - $E_n$-monoidal categories (Definition II.10.7.2),
  - Laplaza $E_n$-monoidal categories (Question A.2.6), and
  - lax $n$-fold monoidal V-categories (Question A.3.1).

As a special case of the last item, the Drinfeld center of a monoidal V-category with V *strict* is studied in [**KYZZ∞, KZ18**]. ◇

**Question A.3.3** (Autonomous Enriched Monoidal Categories)**.** The definition of a monoidal V-category K in [**MP19**, 2.1] assumes that V is braided strict monoidal and K is strict monoidal, so it is the special case of Definition 1.4.2 with both V and K strict. The main theorem in [**MP19**, 1.1] shows that, using their definition, there is a bijective correspondence between (i) some autonomous monoidal V-categories and (ii) some braided oplax monoidal functors from V to the Drinfeld center of an autonomous monoidal category. This bijective correspondence is extended to a 2-equivalence between 2-categories in [**Del∞**, 1.2]. In bicategorical language (Definition I.6.3.9), *autonomy* means that each object $X \in \mathsf{K}$ is equipped with both a left adjoint $X_*$ and a right adjoint $X^*$ that satisfy the triangle identities. A (braided) monoidal functor $(F, F^2, F^0)$ is *oplax* if its monoidal constraint $F^2$ and unit constraint $F^0$ go in the opposite directions as those in Definitions 1.1.6 and 1.1.17, with appropriately adjusted axioms. Discussion of autonomous monoidal categories can be found in [**FY92, JS91b, JS93**].

- Extend the 2-equivalence in [**Del∞**, 1.2] to general monoidal V-categories K, with V and K not necessarily strict, as in Definition 1.4.2.

To extend this 2-equivalence to the general nonstrict case, the coherence and strictification results of enriched monoidal categories in Sections 2.5 and 2.6 will likely be necessary. ◇

## A.4. Homotopy Theory

**Question A.4.1** (Homotopy Theory of Matrix Bicategories)**.** In Example I.8.15.5, we listed some examples of tight symmetric bimonoidal categories C, to which the Bicategorification Theorem I.8.15.4 may be applied to yield a symmetric monoidal bicategory $\mathsf{Mat}^{\mathsf{C}}$.

- What can be said about the homotopy theoretic properties of any of these symmetric monoidal bicategories?
- Consider the previous question for the symmetric bimonoidal bicategories in Question A.1.2.

For instance, for the finite ordinal category $\Sigma$ in Section I.2.4, $\mathsf{Mat}^{\Sigma}$ may be related to a remark in [**JO12**] about the multiplicative structure on the categorical model for the sphere spectrum. ◇

**Question A.4.2** (Categorical Model of $BP$)**.** The Brown-Peterson spectrum $BP$ has an $E_4$ structure [**BM13**], but not an $E_\infty$ structure at any prime [**Law18, Sen∞**]. By Corollary 13.5.2, the Elmendorf-Mandell $K$-theory of each small $E_n$-monoidal category is an $E_n$-symmetric spectrum.

- Is there a small $E_4$-monoidal category (Definition II.10.7.2) whose $K$-theory is $BP$?

A positive answer to this question would provide a categorical model of $BP$. If Question A.5.7 has a positive answer for the 4-fold monoidal category operad $\mathsf{Mon}^4$, then this question also has a positive answer, at least up to weak equivalences. ⋄

**Question A.4.3** (Boardman-Vogt $E_n$-Operads)**.** In a symmetric monoidal category $\mathsf{C}$, the commutative operad $\mathsf{Com}$ is the operad with each entry the monoidal unit $\mathbb{1}$ and structure morphisms given by the coherence isomorphism $\mathbb{1} \otimes \mathbb{1} \cong \mathbb{1}$. The algebras of $\mathsf{Com}$ are precisely commutative monoids. One model of an $E_\infty$-operad (Definition 11.6.1) is $\mathsf{WCom}$, where $\mathsf{W}$ is the Boardman-Vogt $\mathsf{W}$-construction. In the topological setting, it was introduced in [**BV73**]. For a conceptual presentation of the $\mathsf{W}$-construction in terms of coends in a general symmetric monoidal category, see [**Yau20**, Ch. 6–7].

- Describe $E_n$-operads (Definition 12.2.3) for $n \geq 1$ as a filtration of sub-operads of $\mathsf{WCom}$.
- Compare these models of $E_n$-operads to
  - the $n$-fold monoidal category operad $\mathsf{Mon}^n$ (Definition 13.1.12),
  - Berger's filtration [**Ber96**] of the Barratt-Eccles operad $N(E\mathsf{As})$ (Definition 11.4.10),
  - Smith's filtration [**Smi89**] of the Barratt-Eccles operad,
  - Batanin's $E_n$-operads [**Bat07, Bat08**],
  - Fiedorowicz's $E_n$-operads [**Fie∞b**], and
  - the Fulton-MacPherson $E_n$-operads [**Fre17, FM94, GJ∞, Sal01**].

In [**Yau20**, 3.2.11 and 6.3.1], the $\mathsf{W}$-construction of an operad in a symmetric monoidal category is defined entrywise as a coend indexed by a substitution category whose objects are trees and whose morphisms correspond to tree substitution. This applies, in particular, to the commutative operad $\mathsf{Com}$. So describing $E_n$-operads as sub-operads of $\mathsf{WCom}$ would provide a combinatorial description of $E_n$-operads in terms of trees. ⋄

## A.5. Algebraic $K$-Theory

**Question A.5.1** (Multifunctorial $K$-Theory of Pointed Multicategories)**.** Contrary to the claim in [**EM09**, Theorem 1.3], Elmendorf-Mandell $J$-theory $\mathsf{J}^{\mathsf{EM}}$ (Definition 10.3.25) does *not* extend to a multifunctor on all of $\mathsf{Multicat}_*$, the category of small pointed multicategories, but only to the full subcategory $\mathsf{Mod}^{\mathcal{M}\underline{1}}$ of left $\mathcal{M}\underline{1}$-modules, via the symmetric monoidal $\mathsf{Cat}_*$-functor $\mathsf{J}^{\mathcal{T}}$ (Theorem 10.3.17). Examples 10.2.8 and 10.2.9 present some small pointed multicategories that are not left $\mathcal{M}\underline{1}$-modules.

- Is there a $K$-theory multifunctor that is objectwise equivalent to Segal $K$-theory $\mathsf{K}^{\mathsf{Se}}$ and extends to $\mathsf{Multicat}_*$ via the endomorphism multicategory $\mathsf{End}$ in Corollary 5.3.9 and Definition 6.5.1?

The key issue is about the monoidal units. In $\mathsf{Multicat}_*$ the smash unit is $\mathsf{S} = \mathsf{I} \coprod \mathsf{T}$ (Definition 5.6.18), which is different from the monoidal unit $\mathcal{M}\underline{1}$ in $\mathsf{Mod}^{\mathcal{M}\underline{1}}$. Unlike Definition 10.3.16 with $\mathsf{J}^{\mathcal{T}}(\mathcal{M}\underline{1})$, the object $\mathsf{J}^{\mathcal{T}}(\mathsf{S})$ is the terminal $\mathcal{G}_*$-category $*$. If the monoidal unit constraint $(\mathsf{J}^{\mathcal{T}})^0$ for $\mathsf{J}^{\mathcal{T}}$ were to be defined as the unique morphism $J \longrightarrow \mathsf{J}^{\mathcal{T}}(\mathsf{S}) = *$ to the terminal object, as stated in the last paragraph in [**EM09**, Section 5], then $\mathsf{J}^{\mathcal{T}}$ cannot satisfy the unity axioms (1.1.10). The reason is that a general left or right unit isomorphism for $\mathsf{J}^{\mathcal{T}}(-)$ does not factor through the zero morphism in $\mathcal{G}_*$-Cat. So with $J \longrightarrow *$ as the unit constraint, $\mathsf{J}^{\mathcal{T}}$ would not be a monoidal functor. ◇

**Question A.5.2** (Comparison of $\mathsf{K}^{\mathsf{Se}}$ and $\mathsf{K}^{\mathsf{EM}}$ for $\mathcal{M}\underline{1}$-modules)**.** The Segal and Elmendorf-Mandell $K$-theory constructions are defined as the following composites, respectively:

$$\mathsf{K}^{\mathsf{Se}} = \mathsf{K}^{\mathcal{F}} N_* \mathsf{J}^{\mathcal{M}} \mathsf{End} \quad \text{and} \quad \mathsf{K}^{\mathsf{EM}} = \mathsf{K}^{\mathcal{G}} N_* \mathsf{J}^{\mathcal{T}} \mathsf{End}.$$

The domain of $\mathsf{J}^{\mathcal{M}}$ is the category of small pointed multicategories, $\mathsf{Multicat}_*$, and the domain of $\mathsf{J}^{\mathcal{T}}$ is the category of left $\mathcal{M}\underline{1}$-modules within $\mathsf{Multicat}_*$. Thus both $\mathsf{K}^{\mathsf{Se}}$ and $\mathsf{K}^{\mathsf{EM}}$ can be expanded to $\mathsf{Mod}^{\mathcal{M}\underline{1}}$. We will write

$$\mathsf{Mod}^{\mathcal{M}\underline{1}} \xrightarrow{\tilde{\mathsf{K}}^{\mathsf{Se}} = \mathsf{K}^{\mathcal{F}} N_* \mathsf{J}^{\mathcal{M}}} \mathsf{SymSp} \quad \text{and} \quad \mathsf{Mod}^{\mathcal{M}\underline{1}} \xrightarrow{\tilde{\mathsf{K}}^{\mathsf{EM}} = \mathsf{K}^{\mathcal{G}} N_* \mathsf{J}^{\mathcal{T}}} \mathsf{SymSp}.$$

- Is there a (natural) level equivalence $\tilde{\mathsf{K}}^{\mathsf{Se}}\mathsf{P} \longrightarrow \tilde{\mathsf{K}}^{\mathsf{EM}}\mathsf{P}$ for each left $\mathcal{M}\underline{1}$-module $\mathsf{P}$?

The level equivalence $\mathsf{K}^{\mathsf{Se}}\mathsf{C} \longrightarrow \mathsf{K}^{\mathsf{EM}}\mathsf{C}$ given in Theorem 10.6.10 for each small permutative category $\mathsf{C}$ depends crucially on the adjunctions of Proposition 10.6.7 and these, in turn, depend on Proposition 8.5.4, which gives a strong symmetric monoidal adjunction

$$\mathsf{Cat}_*(a, \mathsf{C}) \underset{R}{\overset{L}{\rightleftarrows}} \mathsf{Multicat}_*(\mathcal{M}a, \mathsf{End}(\mathsf{C}))$$

for each small permutative category $\mathsf{C}$ and each pointed finite set $a$. Therefore the proof of Theorem 10.6.10 does not immediately generalize to $\tilde{\mathsf{K}}^{\mathsf{Se}}$ and $\tilde{\mathsf{K}}^{\mathsf{EM}}$. ◇

**Question A.5.3** ($K$-Theory of Matrix Bicategories)**.** There is a $K$-theory construction via the following composite on objects only.

$$\mathsf{Bi}^{\mathsf{tsy}}_{\mathsf{i}} \xrightarrow{\mathsf{Mat}} \mathsf{SMB} \xrightarrow{\text{strictify}} \mathsf{PGray} \xrightarrow{K} \Gamma 2\mathsf{Cat}$$

$\mathsf{Mat}$ is the matrix construction in Theorem I.8.15.4 that sends a tight symmetric bimonoidal category to a symmetric monoidal bicategory. The middle arrow is the strictification of symmetric monoidal bicategories to permutative Gray monoids. The right arrow is the $K$-theory of permutative Gray monoids in [**GJO17b**].

- How does this compare with the $K$-theory of tight symmetric bimonoidal categories in [**BDR04**], which is defined using the direct sum instead of the tensor product of matrices?
- What extra structures on $\Gamma$-2-categories are there when it is the $K$-theory of the strictification of $\mathsf{Mat}^{\mathsf{C}}$ for some tight symmetric bimonoidal category $\mathsf{C}$, such as those in Example I.8.15.5 and $\mathsf{Vect}^{\mathbb{C}}_{\mathsf{c}}$ in Example I.2.5.9? ◇

**Question A.5.4** ($K$-Theory of Matrix Symmetric Bimonoidal Bicategories)**.** Repeat Question A.5.3 for

- the matrix symmetric bimonoidal bicategories in Question A.1.2 and
- the bipermutative Gray monoids in Question A.1.4.

◇

**Question A.5.5** ($K$-Theory of Matrix Permutative Gray Monoids)**.** There is also a $K$-theory construction

$$\mathsf{Bi}^{\mathsf{stsy}} \xrightarrow{\mathsf{Mat}} \mathsf{PGray} \xrightarrow{K} \Gamma\mathsf{2Cat}$$

from *strict* symmetric bimonoidal categories as in Definition I.9.1.1 that uses the matrix construction in Theorem I.9.3.16.

- Can these $K$-theories detect the same weak homotopy types as the ones for tight symmetric bimonoidal categories in Question A.5.3?

This is related to Question A.1.6 in the following sense. It may be the case that strict symmetric bimonoidal categories are too strict to model all categorical equivalence types of tight symmetric bimonoidal categories, leading to a negative answer of Question A.1.6. However, via the matrix construction and $K$-theory, the $\Gamma$-2-categories of strict symmetric bimonoidal categories may model all weak homotopy types of the $\Gamma$-2-categories of tight symmetric bimonoidal categories. ⋄

**Question A.5.6** ($K$-Theory of Distortion Categories)**.** Recall the finite ordinal category $\Sigma$ in Definition I.2.4.1. Here we consider $\Sigma$ as a permutative category with respect to its additive structure $\oplus$. Quillen's +-construction $(B\Sigma)^+$ of the classifying space $B\Sigma$ is the sphere spectrum by the Barratt-Priddy-Quillen Theorem [**BP72**]. A different way to say this is that the algebraic $K$-groups of $\Sigma$ are the stable homotopy groups of the spheres.

- Can the algebraic $K$-groups of
  - the distortion category $\mathcal{D}$ in Section I.4.2,
  - the additive distortion category $\mathcal{D}^{\mathsf{ad}}$ in Section I.4.5, and
  - the braided distortion category $\mathcal{D}^{\mathsf{br}}$ in Section II.5.2

  be computed in similar terms?

By Examples 11.3.18, 11.6.13, and 12.5.4, respectively, $\mathsf{K}^{\mathsf{EM}}\mathcal{D}^{\mathsf{ad}}$, $\mathsf{K}^{\mathsf{EM}}\mathcal{D}$, and $\mathsf{K}^{\mathsf{EM}}\mathcal{D}^{\mathsf{br}}$ are strict ring, $E_\infty$-, and $E_2$-symmetric spectra. Moreover, each of the distortion categories $\mathcal{D}$, $\mathcal{D}^{\mathsf{ad}}$, and $\mathcal{D}^{\mathsf{br}}$ is a Grothendieck construction over $\Sigma$ by, respectively, Propositions I.4.6.5, I.4.6.7, and II.5.5.3.

- Does that yield a computation of their $(B?)^+$ and algebraic $K$-groups? ⋄

**Question A.5.7** (Lifting $K$-Theory Equivalences to Algebras)**.** The Segal $K$-theory functor in Definition 8.5.1 induces an equivalence of homotopy categories via Quillen equivalences, from permutative categories to connective symmetric spectra.

- Do the Quillen equivalences in Segal's $K$-theory lift to the categories of algebras over categorical operads, such as the $E_2$-operad $\mathsf{Br}$ and the $E_n$-operads $\mathsf{Mon}^n$ in Theorems 12.2.4 and 13.2.1?

For a categorical operad $\mathsf{P}$ and a $\mathsf{Cat}$-enriched multicategory $\mathsf{M}$, such as $\mathsf{PermCat}^{\mathsf{su}}$ in Section 6.6, a $\mathsf{P}$-*algebra in* $\mathsf{M}$ is defined as a $\mathsf{Cat}$-enriched multifunctor

$$F : \mathsf{P} \longrightarrow \mathsf{M}.$$

This is equivalent to a $\mathsf{Cat}$-enriched operad morphism

$$\mathsf{P} \longrightarrow \mathsf{End}(A)$$

to the $\mathsf{Cat}$-enriched endomorphism operad of the object $A = F(*)$, with $*$ the unique object in the multicategory $\mathsf{P}$. If the answer to this question is yes for a categorical operad $\mathsf{P}$, then Segal's $K$-theory induces an equivalence between the homotopy categories of (i) $\mathsf{P}$-algebras in permutative categories and (ii) $\mathsf{P}$-algebras in connective symmetric spectra.

To answer this question, it is tempting to use [**WY19**, Th. 4.4 and 4.6], which give sufficient conditions under which a Quillen equivalence between monoidal model categories lifts to a Quillen equivalence between the categories of algebras over some colored operads. This will *not* work because the domain of Segal's $K$-theory is the multicategory $\mathsf{PermCat}^{\mathsf{su}}$ (Section 6.6), which is not a symmetric monoidal category, hence also not a monoidal model category. $\diamond$

# List of Main Facts

**Part I.1**. Symmetric Bimonoidal Categories

**Chapter I.1**. Basic Category Theory

(1.1.11) An adjunction satisfies the triangle identities.
(p. I.11) A functor is an equivalence if and only if it is fully faithful and essentially surjective.
(1.1.14) Left adjoints preserve colimits. Right adjoints preserve limits.
(1.2.1) A monoidal category satisfies the unity axiom and the pentagon axiom.
(1.2.7) A monoidal category satisfies the left and the right unity properties.
(1.3.3) **Mac Lane's Coherence Theorem**. Any two words of the same length in a monoidal category are connected by a unique canonical map.
(1.3.5) **Mac Lane's Strictification Theorem**. Each monoidal category is adjoint equivalent to a strict monoidal category via strong monoidal functors.
(1.3.8) **Symmetric Coherence Theorem**. Any two permuted words of the same length in a symmetric monoidal category are connected by a unique permuted canonical map.
(1.3.10) **Symmetric Strictification Theorem**. Each symmetric monoidal category is adjoint equivalent to a permutative category via strong symmetric monoidal functors.
(1.3.12) **Epstein's Coherence Theorem**. For each (symmetric) monoidal functor $F : \mathsf{C} \longrightarrow \mathsf{D}$ and $F$-iterates $G, H : \mathsf{C}^n \longrightarrow \mathsf{D}$, there exists at most one $F$-coherent map $G \longrightarrow H$.

**Chapter I.2**. Symmetric Bimonoidal Categories

(2.1.2) A symmetric bimonoidal category has two symmetric monoidal structures, left/right multiplicative zero natural isomorphisms, and left/right distributivity natural monomorphisms, and satisfies 24 axioms.
(2.1.2) A bimonoidal category is defined in the same way as a symmetric bimonoidal category, but without the multiplicative symmetry $\xi^{\otimes}$ and the two axioms that involve $\xi^{\otimes}$. So a bimonoidal category is defined by the other 22 axioms.
(2.1.32) There is a tight symmetric bimonoidal category $\mathsf{Vect}^{\mathbb{C}}$ of finite dimensional complex vector spaces.
(2.2.13) Half of the 24 axioms in a symmetric bimonoidal category are formal consequences of the others.
(2.2.14) One axiom is redundant in a bimonoidal category, which is, therefore, determined by 21 axioms.
(2.3.2) Each distributive symmetric monoidal category yields a tight symmetric bimonoidal category, whose sum is the coproduct.

(2.3.3–2.3.5) Symmetric monoidal closed categories with finite coproducts, the category of modules over a commutative ring, and distributive categories are examples of tight symmetric bimonoidal categories.
(2.4.8) The category $\Sigma$ of finite ordinals and permutations is a tight symmetric bimonoidal category.
(2.4.23) The variant $\Sigma'$ of $\Sigma$ is a tight symmetric bimonoidal category.
(2.5.7) Each right bipermutative category is a tight symmetric bimonoidal category.
(2.5.8–2.5.9) $\Sigma'$ and $\mathsf{Vect}_c^{\mathbb{C}}$ are right bipermutative categories.
(2.5.16) Each left bipermutative category is a tight symmetric bimonoidal category.
(2.5.17) $\Sigma$ is a left bipermutative category.
(2.6.2) There is a symmetric bimonoidal groupoid $\Pi$ with syntax of finite types as objects and $\Pi$-terms and $\Pi$-combinators as morphisms.

## **Chapter I.3**. Coherence of Symmetric Bimonoidal Categories

(3.1.6) In the elementary graph $\mathsf{Gr}^{\mathsf{el}}(X)$, $\delta^l$ and $\delta^r$ do not have formal inverses.
(3.1.8) Each prime edge involves at most one nonidentity elementary edge.
(3.1.9) The graph $\mathsf{Gr}(X)$ consists of the vertex set $X^{\mathsf{fr}}$ and prime edges.
(3.1.14) Each functor $\varphi : X \longrightarrow \mathsf{Ob}(\mathsf{C})$ extends additively and multiplicatively to a graph morphism $\varphi : \mathsf{Gr}(X) \longrightarrow \mathsf{C}$.
(3.1.18) The value in $\mathsf{C}$ of a path in $\mathsf{Gr}(X)$ is the composite of the images of its constituent prime edges under $\varphi$.
(3.1.25) An element in $X^{\mathsf{fr}}$ is regular if it has the same support as a formal polynomial whose monomials are distinct in the strict $\{\oplus, \otimes\}$-algebra $X^{\mathsf{st}}$, and whose factors in each monomial are distinct elements in $X$.
(3.1.29) Any two elements in $X^{\mathsf{fr}}$ connected by a path in $\mathsf{Gr}(X)$ have the same support, and one of them is regular if and only if the other one is regular.
(3.2.15) For an element in $X^{\mathsf{fr}}$, the size is equal to the rank if and only if it is a sum, with each summand either in $X$ or a product of two elements in $X$.
(3.3.6) Each element in $X^{\mathsf{fr}}$ has a $0^X$-reduction.
(3.3.11–3.3.12) Any two $0^X$-reductions of an element in $X^{\mathsf{fr}}$ have the same codomain and the same value in a symmetric bimonoidal category.
(3.5.32) Each path in $\mathsf{Gr}(X)$ has a $0^X$-reduction.
(3.5.33) Any two parallel paths in $\mathsf{Gr}(X)$ whose domain has the same support as $0^X$ have the same value in a symmetric bimonoidal category.
(3.6.5) An element in $X^{\mathsf{fr}}$ is a polynomial if and only if it is $\delta$-reduced.
(3.6.9) Each element in $X^{\mathsf{fr}}$ has a $\delta$-reduction.
(3.7.19) Each path in $\mathsf{Gr}(X)$ that does not contain $0^X$ has a $(0^X, \delta)$-reduction.
(3.8.5–3.8.7) Each element in $X^{\mathsf{fr}}$ has a $1^X$-reduction. If the original element is $\delta$-reduced, then all of its $1^X$-reductions have the same codomain and the same value in a symmetric bimonoidal category.
(3.8.14) Each $(0^X, \delta)$-free path in $\mathsf{Gr}(X)$ whose (co)domain is $\delta$-reduced has a $1^X$-reduction.
(3.9.1) **Laplaza's First Coherence Theorem**. In each symmetric bimonoidal category $\mathsf{C}$ satisfying a monomorphism assumption, any two parallel paths in $\mathsf{Gr}(X)$ with a regular domain have the same value in $\mathsf{C}$.
(3.9.9–3.9.10) Theorem I.3.9.1 applies to symmetric bimonoidal categories that are flat, in particular, tight.

(3.10.7) **Bimonoidal Coherence Theorem**. In each bimonoidal category $\mathsf{C}$ satisfying a monomorphism assumption, any two parallel paths in $\mathsf{Gr}^{\mathsf{ns}}(X)$ with a nonsymmetric regular domain have the same value in $\mathsf{C}$.
(3.10.8) Theorem I.3.10.7 applies to flat, in particular, tight, bimonoidal categories.

### Chapter I.4. Coherence of Symmetric Bimonoidal Categories II

(4.2.1) In the distortion category, each object is a finite sequence of nonnegative integers, and each morphism is a finite sequence of permutations.
(4.2.5) The distortion category is a groupoid.
(4.2.12) The additive structure of the distortion category is a permutative category.
(4.2.19) The multiplicative structure of the distortion category is a permutative category.
(4.2.29) The distortion category is a left bipermutative category.
(4.3.1) For a path in $\mathsf{Gr}(X)$, its distortion is defined as its value in the distortion category.
(4.4.3) **Laplaza's Second Coherence Theorem**. In each symmetric bimonoidal category $\mathsf{C}$ satisfying a monomorphism assumption, any two parallel paths in $\mathsf{Gr}(X)$ with the same distortion have the same value in $\mathsf{C}$.
(4.4.5) Theorem I.4.4.3 applies to symmetric bimonoidal categories that are flat, in particular, tight.
(4.5.2) In the additive distortion category, each object is a finite sequence of nonnegative integers, and each morphism is a permutation.
(4.5.6) The additive distortion category is a groupoid and a tight bimonoidal category. It faithfully embeds into the distortion category.
(4.5.7) For a path in $\mathsf{Gr}^{\mathsf{ns}}(X)$, its additive distortion is defined as its value in the additive distortion category.
(4.5.8) **Bimonoidal Coherence Theorem II**. In each bimonoidal category $\mathsf{C}$ satisfying a monomorphism assumption, any two parallel paths in $\mathsf{Gr}^{\mathsf{ns}}(X)$ with the same additive distortion have the same value in $\mathsf{C}$.
(4.5.9) Theorem I.4.5.8 applies to flat, in particular, tight, bimonoidal categories.
(4.6.5) The distortion category $\mathcal{D}$ is isomorphic to the Grothendieck construction $\int_{\Sigma} F$.
(4.6.7) The additive distortion category $\mathcal{D}^{\mathsf{ad}}$ is isomorphic to the Grothendieck construction $\int_{\Sigma} F^{\mathsf{ad}}$.

### Chapter I.5. Strictification of Tight Symmetric Bimonoidal Categories

(5.1.1) A symmetric bimonoidal functor is a functor equipped with two symmetric monoidal functor structures, and satisfies two axioms.
(5.1.10) There is a 1-category $\mathsf{Bi}^{\mathsf{sy}}$ of small symmetric bimonoidal categories and symmetric bimonoidal functors.
(5.1.11) Each symmetric monoidal functor between distributive symmetric monoidal categories induces a symmetric bimonoidal functor.
(5.1.15–5.1.16) $\Sigma$ and $\Sigma'$ are isomorphic via symmetric bimonoidal functors.
(5.3.9) Each tight symmetric bimonoidal category $\mathsf{C}$ has an associated right bipermutative category $\mathsf{A}$, whose objects are formal polynomials in the objects in $\mathsf{C}$.
(5.4.6–5.4.7) **Bipermutative Strictification Theorems**. Each tight symmetric bimonoidal category is adjoint equivalent to a right/left bipermutative category via symmetric bimonoidal functors.

(5.5.1) A bimonoidal functor is a functor equipped with an additive symmetric monoidal functor structure and a multiplicative monoidal functor structure that satisfies four axioms.
(5.5.4) There is a 1-category $\mathsf{Bi}$ of small bimonoidal categories and bimonoidal functors.
(5.5.10) Each tight bimonoidal category $\mathsf{C}$ has an associated right rigid bimonoidal category $\mathsf{A}$, whose objects are formal polynomials in the objects in $\mathsf{C}$.
(5.5.11–5.5.12) **Rigid Strictification Theorems**. Each tight bimonoidal category is adjoint equivalent to a right/left rigid bimonoidal category via bimonoidal functors.

**Part I.2**. Bicategorical Aspects of Symmetric Bimonoidal Categories

**Chapter I.6**. Definitions from Bicategory Theory

(6.1.2) A bicategory has objects, (identity) 1-cells, (identity) 2-cells, vertical and horizontal compositions, an associator, and two unitors, and satisfies the unity axiom and the pentagon axiom.
(6.1.8) A 2-category is a bicategory whose associator and unitors are identities.
(6.1.10) A 2-category can be described by data and axioms.
(6.1.11) A locally small 2-category is precisely a $\mathsf{Cat}$-category.
(6.1.16) A monoidal category is a one-object bicategory.
(6.2.1) A lax functor has an object assignment, local functors, and two laxity constraints, and satisfies the lax associativity axiom and the lax unity axioms.
(6.2.11) There is a 1-category $\mathsf{Bicat}$ with small bicategories as objects and lax functors as morphisms.
(6.2.14) A lax transformation has component 1-cells and natural component 2-cells, and satisfies the lax unity axiom and the lax naturality axiom.
(6.2.26) There is a 2-category $\mathsf{2Cat}$ of small 2-categories, 2-functors, and 2-natural transformations.
(6.3.1) A modification has component 2-cells and satisfies the modification axiom.
(6.3.7) For bicategories $\mathsf{B}$ and $\mathsf{B}'$ with $\mathsf{Ob}(\mathsf{B})$ a set, there is a bicategory $\mathsf{Bicat}(\mathsf{B}, \mathsf{B}')$ with lax functors $\mathsf{B} \longrightarrow \mathsf{B}'$ as objects, lax transformations as 1-cells, and modifications as 2-cells. It is a 2-category if $\mathsf{B}'$ is a 2-category. It contains a full subbicategory $\mathsf{Bicat}^{\mathsf{ps}}(\mathsf{B}, \mathsf{B}')$ with pseudofunctors as objects and strong transformations as 1-cells.
(6.3.9) An adjunction in a bicategory consists of two 1-cells and two 2-cells, and satisfies two triangle identities.
(6.4.1) A monoidal bicategory has a base bicategory, a monoidal composition, a monoidal identity, a monoidal associator, two monoidal unitors, a pentagonator, and three 2-unitors, and satisfies the non-abelian 4-cocycle condition and two normalization axioms.
(6.5.3) A braided monoidal bicategory is a monoidal bicategory equipped with a braiding and two hexagonators, and satisfies four axioms.
(6.5.7) A sylleptic monoidal bicategory is a braided monoidal bicategory equipped with a syllepsis that satisfies two axioms.
(6.5.9) A symmetric monoidal bicategory is a sylleptic monoidal bicategory that satisfies the triple braid axiom.
(6.6.12) The 1-category $\mathsf{2Cat}$ equipped with the Gray tensor product is a symmetric monoidal closed category $\mathsf{Gray}$.

(6.6.13) A Gray monoid is a monoid in Gray.
(6.7.1) A permutative Gray monoid is a Gray monoid equipped with a Gray symmetry that satisfies three axioms.
(6.7.16) A permutative 2-category is a monoid in $(2\mathsf{Cat}, \times)$ that is equipped with a symmetry 2-natural isomorphism and that satisfies the same three axioms as for permutative Gray monoids.

## **Chapter I.7**. Baez's Conjecture

(7.1.2) A bimonoidal natural transformation is a natural transformation that is a monoidal natural transformation for each of the additive structure and the multiplicative structure.
(7.1.7) There is a 2-category $\mathsf{Bi}^{\mathsf{sy}}$ of small symmetric bimonoidal categories, symmetric bimonoidal functors, and bimonoidal natural transformations.
(7.1.9) Each monoidal natural transformation between symmetric monoidal functors between distributive symmetric monoidal categories is a bimonoidal natural transformation.
(7.2.9) For each symmetric bimonoidal category $\mathsf{C}$, there is a strong symmetric monoidal functor $F_\oplus : \Sigma \longrightarrow \mathsf{C}$ between the additive structures.
(7.3.28) For each flat symmetric bimonoidal category $\mathsf{C}$, there is a symmetric monoidal functor $F_\otimes : \Sigma \longrightarrow \mathsf{C}$ between the multiplicative structures.
(7.4.4) For each flat symmetric bimonoidal category $\mathsf{C}$, $F : \Sigma \longrightarrow \mathsf{C}$ is a robust symmetric bimonoidal functor.
(7.5.8) Epstein's Coherence Theorem I.1.3.12 has a bimonoidal analogue.
(7.6.2–7.6.3) For a symmetric bimonoidal category $\mathsf{C}$ and robust symmetric bimonoidal functors $G, H : \Sigma \longrightarrow \mathsf{C}$, there is at most one bimonoidal natural transformation $G \longrightarrow H$, which must be invertible if it exists.
(7.7.9) For each flat symmetric bimonoidal category $\mathsf{C}$ and robust symmetric bimonoidal functor $G : \Sigma \longrightarrow \mathsf{C}$, there exists a unique bimonoidal natural transformation $\theta : F \longrightarrow G$, which is, moreover, invertible.
(7.8.1) **Baez's Conjecture**. $\Sigma$ is a lax bicolimit of the 2-functor $\varnothing \longrightarrow \mathsf{Bi}_{\mathsf{r}}^{\mathsf{fsy}}$.
(7.8.3) **Baez's Conjecture, Ver. 2**. $\Sigma'$ is a lax bicolimit of the 2-functor $\varnothing \longrightarrow \mathsf{Bi}_{\mathsf{r}}^{\mathsf{fsy}}$.

## **Chapter I.8**. Symmetric Monoidal Bicategorification

(8.1.1) For a category $\mathsf{C}$, $\mathsf{Mat}^{\mathsf{C}}_{m,n}$ has $n \times m$ matrices of objects in $\mathsf{C}$ as objects and $n \times m$ matrices of morphisms in $\mathsf{C}$ as morphisms.
(8.1.8) The matrix product is a functor.
(8.2.2) There is a natural isomorphism $\ell_A : \mathbb{1}^n A \xrightarrow{\cong} A$ for each flat bimonoidal category $\mathsf{C}$ and $A \in \mathsf{Mat}^{\mathsf{C}}_{m,n}$.
(8.2.8) There is a natural isomorphism $r_A : A\mathbb{1}^m \xrightarrow{\cong} A$ for each flat bimonoidal category $\mathsf{C}$ and $A \in \mathsf{Mat}^{\mathsf{C}}_{m,n}$.
(8.3.1) There is a natural isomorphism $a : (CB)A \xrightarrow{\cong} C(BA)$ for each tight bimonoidal category $\mathsf{C}$.
(8.4.12) For each tight bimonoidal category $\mathsf{C}$, $\mathsf{Mat}^{\mathsf{C}}$ is a bicategory.
(8.4.14) $\mathsf{Mat}^{\mathsf{C}}_{n,n}$ is a monoidal category.
(8.6.7) The matrix tensor product is a functor.
(8.7.31) The triple $(\boxtimes, \boxtimes^2, \boxtimes^0)$ on $\mathsf{Mat}^{\mathsf{C}}$ is a pseudofunctor.
(8.8.49) The quadruple $(a^{\boxtimes}, a^{\boxtimes\bullet}, \eta^a, \varepsilon^a)$ is an adjoint equivalence.

(8.9.9) The quadruple $(\ell^{\boxtimes}, \ell^{\boxtimes\bullet}, \eta^{\ell}, \varepsilon^{\ell})$ is an adjoint equivalence.
(8.9.21) The quadruple $(r^{\boxtimes}, r^{\boxtimes\bullet}, \eta^{r}, \varepsilon^{r})$ is an adjoint equivalence.
(8.10.4) $\pi$ is an invertible modification.
(8.11.4) $\mu$ is an invertible modification.
(8.11.9) $\lambda^{\boxtimes}$ is an invertible modification.
(8.11.14) $\rho^{\boxtimes}$ is an invertible modification.
(8.12.9) For each tight symmetric bimonoidal category $\mathsf{C}$, $\mathsf{Mat}^{\mathsf{C}}$ is a monoidal bicategory.
(8.13.13) For $\sigma \in \Sigma_m$, there is a natural isomorphism $r^{\sigma}_{A} : A\mathbb{1}^{\sigma} \xrightarrow{\cong} A^{\sigma}$ for $A \in \mathsf{Mat}^{\mathsf{C}}_{m,n}$.
(8.13.16) For $\theta \in \Sigma_n$, there is a natural isomorphism $\ell^{\theta}_{A} : \mathbb{1}^{\theta}A \xrightarrow{\cong} {}_{\theta^{-1}}A$ for $A \in \mathsf{Mat}^{\mathsf{C}}_{m,n}$.
(8.13.20) The matrix tensor products $A \boxtimes B$ and $B \boxtimes A$ differ by a column permutation, a row permutation, and a multiplicative symmetry in each entry.
(8.13.44) The quadruple $(\beta, \beta^{\bullet}, \eta^{\beta}, \varepsilon^{\beta})$ is an adjoint equivalence.
(8.14.12) $R_{-|--}$ is an invertible modification.
(8.14.24) $R_{--|-}$ is an invertible modification.
(8.14.26) For each tight symmetric bimonoidal category $\mathsf{C}$, $\mathsf{Mat}^{\mathsf{C}}$ is a braided monoidal bicategory.
(8.15.4) **Bicategorification Theorem**. For each tight symmetric bimonoidal category $\mathsf{C}$, $\mathsf{Mat}^{\mathsf{C}}$ is a symmetric monoidal bicategory.
(8.15.5) Coordinatized 2-vector spaces $2\mathsf{Vect}_{\mathsf{c}} = \mathsf{Mat}^{\mathsf{Vect}^{\mathbb{C}}}$ form a symmetric monoidal bicategory.

**Chapter I.9**. Matrix Permutative Gray Monoids

(9.1.1) A symmetric bimonoidal category is strict if (i) the additive and the multiplicative structures are both permutative categories, and (ii) $\xi^{\oplus}$, $\delta^{l}$, $\delta^{r}$, $\lambda^{\bullet}$, and $\rho^{\bullet}$ are identities.
(9.1.7–9.1.8) $\mathsf{Mat}^{\mathsf{C}}$ is a 2-category if $\mathsf{C}$ is a strict (symmetric) bimonoidal category.
(9.2.12) For each strict symmetric bimonoidal category $\mathsf{C}$, $\mathsf{Mat}^{\mathsf{C}}$ is a Gray monoid.
(9.3.16) For each strict symmetric bimonoidal category $\mathsf{C}$, $\mathsf{Mat}^{\mathsf{C}}$ is a permutative Gray monoid.
(9.4.2) For each strict symmetric bimonoidal category $\mathsf{C}$ with $\xi^{\otimes}$ the identity, $\mathsf{Mat}^{\mathsf{C}}$ is a permutative 2-category.

**Part II.1**. Braided Bimonoidal Categories

**Chapter II.1**. Preliminaries on Braided Structures

(1.1.1) The braid group $B_n$ on $n$ strings is generated by $s_1, \ldots, s_{n-1}$ and subject to two braid relations.
(1.1.9) Sum braids generalize block sums of permutations.
(1.1.12) Each braid has an underlying permutation.
(1.1.20) Block braids generalize block permutations.
(1.2.4) Elementary block braids generalize interval-swapping permutations.
(1.2.14) Elementary block braids are compatible with sum braids.
(1.2.16) Elementary block braids satisfy the hexagon axioms.
(1.3.15) A braided monoidal category is a monoidal category equipped with a braiding that satisfies two hexagon axioms.

(1.3.18) A braided monoidal functor is defined in the same way as a symmetric monoidal functor.
(1.3.21) In each braided monoidal category, the left unit isomorphism uniquely determines the right unit isomorphism, and vice versa, via the braiding.
(1.3.28, 1.3.31) Each braided monoidal category satisfies the third Reidemeister move.
(1.3.36) A symmetric monoidal category is precisely a braided monoidal category whose braiding satisfies the symmetry axiom.
(1.4.27) The Drinfeld center of a monoidal category is a braided monoidal category.
(1.5.3) The symmetric center of a braided monoidal category is a symmetric monoidal category.
(1.6.3) **Braided Coherence Theorem**. Two braided canonical maps with the same (co)domain in a braided monoidal category are equal if their underlying braids are equal.
(1.6.5) **Braided Strictification Theorem**. Each braided monoidal category is adjoint equivalent to a braided strict monoidal category via strong braided monoidal functors.

## **Chapter** II.2. Braided Bimonoidal Categories

(2.1.29) A braided bimonoidal category is a category equipped with an additive symmetric monoidal structure, a multiplicative braided monoidal structure, left/right multiplicative zero natural isomorphisms, and left/right distributivity natural monomorphisms, and satisfies twelve Laplaza's axioms and two additional axioms involving the braiding.
(2.1.37) Tight braided bimonoidal categories are equivalent to BD categories in the sense of Blass and Gurevich.
(2.2.1) Each braided bimonoidal category satisfies all 24 Laplaza axioms.
(2.2.3) A symmetric bimonoidal category is precisely a braided bimonoidal category whose braiding satisfies the symmetry axiom.
(2.3.2) In an Ab-category, composition with a zero morphism yields a zero morphism, and composition commutes with taking the additive inverse of a morphism.
(2.3.7 (1)) For any two objects in an Ab-category, a product, a coproduct, and a direct sum are equivalent.
(2.3.7 (2)) In an Ab-category, the direct sum morphism $f \oplus f'$ can be characterized in terms of the inclusions.
(2.3.7 (3)) In an Ab-category, the sum morphism $f + g$ factors as $\nabla_B(f \oplus g)\Delta_A$.
(2.3.7 (4)) A functor between Ab-categories whose domain has all direct sums is an additive functor if and only if it preserves direct sums.
(2.3.12) For any two objects in a preadditive category, the zero morphism is the unique morphism that factors through the zero object.
(2.3.15) An abelian category is an Ab-category with a zero object, a direct sum for any two objects, and a (co)kernel for each morphism, such that (i) each monomorphism is a kernel and (ii) each epimorphism is a cokernel.
(2.3.17) Each abelian category has all finite (co)limits, with (co)products given by direct sums.
(2.4.22) An abelian category with a compatible braided monoidal structure is a tight braided bimonoidal category.

(2.5.1–2.5.2) An abelian category with a compatible (symmetric) monoidal structure is a tight (symmetric) bimonoidal category.

**Chapter II.3**. Applications to Quantum Groups and Topological Quantum Computation

(3.1.19) A braided bialgebra is a bialgebra equipped with an $R$-matrix that satisfies two axioms. A symmetric bialgebra is a braided bialgebra in which the inverse of the $R$-matrix is its opposite.
(3.1.27–3.1.30) Each cocommutative bialgebra is a symmetric bialgebra with the $R$-matrix $1 \otimes 1$. Examples include group bialgebras, the universal enveloping bialgebra of a Lie algebra, and Sweedler's 4-dimensional non-(co)commutative bialgebra.
(3.1.33) Each anyonic quantum group is a braided bialgebra.
(3.2.6) The category of left modules over each bialgebra is a monoidal category under the tensor product.
(3.2.12) The category of left modules over each braided bialgebra is a braided monoidal category.
(3.2.13) The category of left modules over each symmetric bialgebra is a symmetric monoidal category.
(3.2.19) The category of left modules over each bialgebra is a tight bimonoidal category. The braided and the symmetric analogues are also true.
(3.3.27) The Fibonacci anyons form a monoidal category.
(3.4.5) The Fibonacci anyons form a braided monoidal category.
(3.4.13) The Fibonacci anyons form a tight braided bimonoidal category.
(3.5.27) The Ising anyons form a monoidal category.
(3.6.7) The Ising anyons form a braided monoidal category.
(3.6.14) The Ising anyons form a tight braided bimonoidal category.

**Chapter II.4**. Bimonoidal Centers

(4.2.6) The additive structure of the bimonoidal Drinfeld center is a symmetric monoidal category.
(4.3.3) The multiplicative structure of the bimonoidal Drinfeld center is a braided monoidal category.
(4.4.3) For each tight bimonoidal category, the bimonoidal Drinfeld center is a tight braided bimonoidal category.
(4.5.3) For each braided bimonoidal category, the bimonoidal symmetric center is a symmetric bimonoidal category.

**Chapter II.5**. Coherence of Braided Bimonoidal Categories

(5.1.2) A left permbraided category has an additive permutative structure, a multiplicative braided strict monoidal structure, and identities for $\lambda^{\bullet}$, $\rho^{\bullet}$, $\delta^{l}$, $\xi^{\otimes}_{-,0}$, and $\xi^{\otimes}_{0,-}$, and satisfies four braided bimonoidal category axioms.
(5.1.8) Each left bipermutative category is a left permbraided category.
(5.1.10) Each left permbraided category is a tight braided bimonoidal category.
(5.1.11) A right permbraided category has an additive permutative structure, a multiplicative braided strict monoidal structure, and identities for $\lambda^{\bullet}$, $\rho^{\bullet}$, $\delta^{r}$, $\xi^{\otimes}_{-,0}$, and $\xi^{\otimes}_{0,-}$, and satisfies four braided bimonoidal category axioms.
(5.1.17) Each right bipermutative category is a right permbraided category.
(5.1.19) Each right permbraided category is a tight braided bimonoidal category.

(5.2.7) The braided distortion category is a groupoid.
(5.2.13) The additive structure of the braided distortion category is a permutative category.
(5.2.21) The multiplicative structure of the braided distortion category is a braided strict monoidal category.
(5.2.28) In the braided distortion category, the right distributivity morphism $\delta^r$ has identity braid components.
(5.2.30) The braided distortion category is a left permbraided category.
(5.2.33–5.2.34) The braided distortion category is a tight braided bimonoidal category and satisfies all 24 Laplaza axioms.
(5.3.14) For a path in $\mathsf{Gr}(X)$, its value in a braided bimonoidal category is the composite of the images of its constituent prime edges.
(5.3.15) The braided distortion of a path in $\mathsf{Gr}(X)$ is its value in the braided distortion category.
(5.4.4) **Braided Bimonoidal Coherence Theorem**. In each braided bimonoidal category $\mathsf{C}$ satisfying a monomorphism assumption, any two parallel paths in $\mathsf{Gr}(X)$ with the same braided distortion have the same value in $\mathsf{C}$.
(5.4.6) Theorem II.5.4.4 applies to flat, in particular, tight, braided bimonoidal categories.
(5.5.3) The braided distortion category $\mathcal{D}^{\mathsf{br}}$ is isomorphic to the Grothendieck construction $\int_{\Sigma} F^{\mathsf{br}}$.

**Chapter II.6**. Strictification of Tight Braided Bimonoidal Categories

(6.1.1) A braided bimonoidal functor is a functor equipped with an additive symmetric monoidal structure and a multiplicative braided monoidal structure, and satisfies two axioms.
(6.1.10) There is a category $\mathsf{Bi}^{\mathsf{br}}$ with small braided bimonoidal categories as objects and braided bimonoidal functors as morphisms.
(6.1.12) Each braided monoidal functor that is also an additive functor between abelian categories with a compatible braided monoidal structure canonically extends to a braided bimonoidal functor.
(6.1.15) Each symmetric monoidal functor that is also an additive functor between abelian categories with a compatible symmetric monoidal structure canonically extends to a symmetric bimonoidal functor.
(6.2.39) Each tight braided bimonoidal category has a canonically associated right permbraided category.
(6.3.6–6.3.7) **Permbraided Strictification**. Each tight braided bimonoidal category is adjoint equivalent to a right/left permbraided category via braided bimonoidal functors.

**Chapter II.7**. The Braided Baez Conjecture

(7.1.4) There is a 2-category $\mathsf{Bi}^{\mathsf{br}}$ of small braided bimonoidal categories, braided bimonoidal functors, and bimonoidal natural transformations.
(7.1.7) Each monoidal natural transformation between braided monoidal functors that are also additive functors between abelian categories with a compatible braided monoidal structure is a bimonoidal natural transformation.
(7.2.4) For each braided bimonoidal category $\mathsf{C}$, there is a strong symmetric monoidal functor $F_{\oplus} : \Sigma \longrightarrow \mathsf{C}$ between the additive structures.

(7.2.9) For each flat braided bimonoidal category $\mathsf{C}$, there is a braided monoidal functor $F_{\otimes} : \Sigma \longrightarrow \mathsf{C}$ between the multiplicative structures.
(7.2.11) For each flat braided bimonoidal category $\mathsf{C}$, $F : \Sigma \longrightarrow \mathsf{C}$ is a robust braided bimonoidal functor.
(7.3.4) **Braided Baez Conjecture**. $\Sigma$ is a lax bicolimit of the 2-functor $\varnothing \longrightarrow \mathsf{Bi}_{\mathsf{r}}^{\mathsf{fbr}}$.
(7.3.6) **Braided Baez Conjecture, Ver. 2**. $\Sigma'$ is a lax bicolimit of the 2-functor $\varnothing \longrightarrow \mathsf{Bi}_{\mathsf{r}}^{\mathsf{fbr}}$.

## Chapter II.8. Monoidal Bicategorification

(8.1.13) For each tight braided bimonoidal category $\mathsf{C}$, $\mathsf{Mat}^{\mathsf{C}}$ is a bicategory.
(8.4.7) For each tight braided bimonoidal category $\mathsf{C}$, $\mathsf{Mat}^{\mathsf{C}}$ is a monoidal bicategory.

# Part II.2. $E_n$-Monoidal Categories

## Chapter II.9. Ring, Bipermutative, and Braided Ring Categories

(9.1.15) Tight ring categories form a subclass of tight bimonoidal categories.
(9.1.19) Right and left rigid bimonoidal categories are tight ring categories.
(9.1.20) Each tight ring category is adjoint equivalent to a right, respectively, left, rigid bimonoidal category.
(9.2.14) Each small permutative category has an endomorphism ring category.
(9.2.20) Each small permutative category has an endomorphism tight ring category.
(9.3.7) Each tight bipermutative category yields a tight symmetric bimonoidal category.
(9.3.12) Right and left bipermutative categories are tight bipermutative categories.
(9.3.13) Each tight bipermutative category is adjoint equivalent to a right, respectively, left, bipermutative category.
(9.4.7) In a bipermutative category, about half of the ring category axioms are redundant.
(9.5.4) A bipermutative category is a braided ring category whose braiding satisfies the symmetry axiom.
(9.5.5) In a braided ring category, about half of the ring category axioms are redundant.
(9.5.6) Tight braided ring categories form a subclass of tight braided bimonoidal categories.
(9.5.10) Right and left permbraided categories are tight braided ring categories.
(9.5.11) Each tight braided ring category is adjoint equivalent to a right, respectively, left, permbraided category.
(9.6.1) The bimonoidal Drinfeld center of a tight ring category is a tight braided ring category.
(9.6.4) The symmetric center of a braided ring category with left factorization a natural epimorphism is a bipermutative category.

## Chapter II.10. Iterated and $E_n$-Monoidal Categories

(10.1.9) A 1-fold monoidal category is a strict monoidal category.
(10.1.14) Braided strict monoidal categories form a subclass of 2-fold monoidal categories.

(10.1.21) Permutative categories form a subclass of $n$-fold monoidal categories for each $n \geq 2$.
(10.2.3) A totally ordered set with a least element forms a permutative category with identity symmetry.
(10.2.8) A totally ordered monoid whose unit is also the least element forms a 2-fold monoidal category.
(10.3.7) A 1-fold monoidal functor is a strictly unital monoidal functor.
(10.3.11) A braided strictly unital monoidal functor is a 2-fold monoidal functor.
(10.3.15) A symmetric strictly unital monoidal functor is an $n$-fold monoidal functor for each $n \geq 2$.
(10.3.20) The composite of two $n$-fold monoidal functors is an $n$-fold monoidal functor.
(10.4.2) $\mathsf{MCat}^n$ is a monoidal category.
(10.4.5) Monoids in $\mathsf{MCat}^n$ are small $(n+1)$-fold monoidal categories.
(10.4.13) A morphism of monoids in $\mathsf{MCat}^n$ is an $(n+1)$-fold monoidal functor with the last monoidal constraint the identity.
(10.5.9) $\mathsf{FMon}^n : \mathsf{Cat} \longrightarrow \mathsf{MCat}^n_{\mathsf{st}}$ is the left adjoint of the forgetful functor.
(10.5.18) $\mathsf{FMon}^n(\mathsf{C})$ decomposes into a coproduct $\coprod_{k\geq 0} \mathsf{Mon}^n(k) \times_{\Sigma_k} \mathsf{C}^{\times k}$.
(10.5.26) $\coprod_{k\geq 0} \mathsf{Mon}^n(k)/\Sigma_k$ is the free $n$-fold monoidal category on one object.
(10.5.28) There are evaluation functors $\theta_k : \mathsf{Mon}^n(k) \times_{\Sigma_k} \mathsf{C}^{\times k} \longrightarrow \mathsf{C}$ for each small $n$-fold monoidal category $\mathsf{C}$.
(10.6.8 (1)) Each morphism set in $\mathsf{Mon}^n(k)$ has at most one morphism.
(10.6.8 (2)) There exists a morphism $A \longrightarrow B \in \mathsf{Mon}^n$ if and only if $a \otimes_i b \in A$ implies either $a \otimes_j b \in B$ for some $j \geq i$ or $b \otimes_j a \in B$ for some $j > i$.
(10.6.9) In each $n$-fold monoidal category, each formal diagram built from identity morphisms, the exchanges $\{\eta^{i,j}\}_{i<j}$, the monoidal products $\{\otimes_i\}_{i=1}^n$, and composites is commutative.
(10.7.13) An $E_1$-monoidal category is a ring category.
(10.8.1) Braided ring categories form a subclass of $E_2$-monoidal categories.
(10.9.1) Bipermutative categories form a subclass of $E_n$-monoidal categories for each $n \geq 2$.
(10.10.2) Each small category has a free $E_n$-monoidal category.

## Part III.1. Enriched Monoidal Categories and Multicategories

### Chapter III.1. Enriched Monoidal Categories

(1.1.31) **Mac Lane's Coherence Theorem**. Any two words of the same length in a monoidal category are connected by a unique canonical map.
(1.1.32) **Mac Lane's Strictification Theorem**. Each monoidal category is adjoint equivalent to a strict monoidal category via strong monoidal functors.
(1.1.38) **Braided Coherence Theorem**. Two braided canonical maps with the same (co)domain in a braided monoidal category are equal if their underlying braids are equal.
(1.1.39) **Braided Strictification Theorem**. Each braided monoidal category is adjoint equivalent to a braided strict monoidal category via strong braided monoidal functors.
(1.1.41) **Symmetric Coherence Theorem**. Any two permuted words of the same length in a symmetric monoidal category are connected by a unique permuted canonical map.

(1.1.42) **Symmetric Strictification Theorem**. Each symmetric monoidal category is adjoint equivalent to a permutative category via strong symmetric monoidal functors.
(1.1.44) **Epstein's Coherence Theorem**. For each (symmetric) monoidal functor $F : \mathsf{C} \longrightarrow \mathsf{D}$ and $F$-iterates $G, H : \mathsf{C}^n \longrightarrow \mathsf{D}$, there exists at most one $F$-coherent map $G \longrightarrow H$.
(1.2.1) For a monoidal category V, a V-category has hom objects in V satisfying associativity and unity axioms.
(1.2.4) A V-functor satisfies composition and identity axioms.
(1.2.7) A V-natural transformation satisfies a naturality axiom.
(1.2.13) There is a 2-category formed by small V-categories, V-functors, and V-natural transformations.
(1.2.16) The opposite of a V-category is defined if V is braided monoidal.
(1.3.3) The tensor product of V-categories is defined if V is braided monoidal.
(1.3.6) The tensor product of V-categories is 2-functorial.
(1.3.35) The underlying 1-category of V-Cat is monoidal if V is braided, and is symmetric monoidal if V is symmetric.
(1.4.2) A monoidal V-category has associator and unitor V-natural transformations that satisfy unity axioms and a pentagon axiom.
(1.4.7) Composition in a monoidal V-category has an enriched interchange.
(1.4.10) The definition of braided monoidal V-category requires that V be symmetric monoidal.
(1.4.10) A braided monoidal V-category has a V-natural braiding satisfying two hexagon axioms.
(1.4.13) A symmetric monoidal V-category is a braided monoidal V-category satisfying an additional symmetry axiom.
(1.4.17) A monoidal V-functor satisfies associativity and unity axioms.
(1.4.18) A braided monoidal V-functor satisfies a braid axiom. A symmetric monoidal V-functor is a braided monoidal V-functor whose domain and codomain are symmetric monoidal V-categories.
(1.4.22) A monoidal V-natural transformation satisfies monoidal naturality and monoidal unity axioms.
(1.4.25) There are 2-categories formed by each of: monoidal V-categories, braided monoidal V-categories, and symmetric monoidal V-categories with, in each case, the corresponding V-functors and V-natural transformations.
(1.5.1) We use the term Cat-monoidal 2-category to indicate monoidal V-categories when V = Cat.
(1.5.2) The underlying 1-category of a plain/braided/symmetric Cat-monoidal 2-category has the corresponding structure as a 1-category.
(1.5.3) A Cat-monoidal 2-category has a strict form of the data and axioms for a monoidal bicategory. Similar statements hold for the braided and symmetric cases.
(1.5.4) With the Cartesian product, Cat is a symmetric Cat-monoidal 2-category.
(1.5.5) For a braided monoidal category V, V-Cat is a Cat-monoidal 2-category. If V is symmetric, then so is V-Cat.

## **Chapter III.2**. Change of Enrichment

(2.1.2) Change of enrichment along a monoidal functor $U$ is 2-functorial.

(2.1.7) The functor from V-Cat to Cat that takes underlying categories is injective on 2-cells.
(2.2.7) The assignment $E : \mathsf{V} \longmapsto \mathsf{V}\text{-Cat}$ is 2-functorial with respect to monoidal functors and monoidal natural transformations.
(2.2.11) For small monoidal V, there is a 2-equivalence between V-Cat and $(\mathsf{V}_{\mathsf{st}})$-Cat.
(2.3.7) If $U$ is braided moniodal, then the change of enrichment induced by $U$ is a Cat-monoidal 2-functor. If $U$ is symmetric, then so is the change of enrichment.
(2.3.9) The assignment $E : \mathsf{V} \longmapsto \mathsf{V}\text{-Cat}$ of braided monoidal categories to Cat-monoidal 2-categories is 2-functorial. A similar result holds for symmetric monoidal categories with $E$ producing symmetric Cat-monoidal categories.
(2.3.16) For small symmetric monoidal V, there is a symmetric Cat-monoidal 2-equivalence between V-Cat and $(\mathsf{V}_{\mathsf{st}})$-Cat.
(2.4.10, 2.4.15) For braided monoidal $U : \mathsf{V} \longrightarrow \mathsf{W}$, change of enrichment along $U$ induces 2-functors between the 2-categories of monoidal V- and W-categories. If $U$ is symmetric, then a similar result holds for braided and symmetric monoidal V- and W-categories.
(2.4.17) The underlying category of a monoidal V-category is monoidal. Similar statements hold for braided and symmetric cases, and for functors and natural transformations.
(2.5.1) Given the data of a monoidal V-category, the enriched monoidal category axioms are satisfied if and only if the underlying data satisfy the ordinary monoidal category axioms. Similar results hold for the braided and symmetric monoidal cases, and also for functors and natural transformations.
(2.5.6) **Enriched (Braided/Symmetric) Monoidal Coherence Theorem**. Any two V-words of the same length in a monoidal V-category are connected by a unique canonical V-map. Similar coherence results hold in the braided and symmetric cases.
(2.5.8) **Enriched Epstein's Coherence Theorem**. For each (symmetric) monoidal V-functor $F : \mathsf{K} \longrightarrow \mathsf{L}$, and $F$-iterates $G, H : \mathsf{K}^{\otimes n} \longrightarrow \mathsf{L}$, there exists at most one $F$-coherent map $G \longrightarrow H$.
(2.6.1) **Enriched Monoidal Strictification Theorem**. Each monoidal V-category is adjoint V-equivalent to a strict monoidal V-category via strong monoidal V-functors.
(2.6.3) **Enriched Braided Strictification Theorem**. Each braided monoidal V-category is adjoint V-equivalent to a braided strict monoidal V-category via strong braided monoidal V-functors.
(2.6.4) **Enriched Symmetric Strictification Theorem**. Each symmetric monoidal V-category is adjoint V-equivalent to a strict monoidal V-category via strong symmetric monoidal V-functors.

## Chapter III.3. Self-Enrichment and Enriched Yoneda

(3.1.11) Each symmetric monoidal closed V has a canonical enrichment over itself, $\underline{\mathsf{V}}$.
(3.2.1, 3.2.2) A category enriched over symmetric monoidal closed V has co/represented V-functors $\mathcal{Y}^X$ and $\mathcal{Y}_Y$ to $\underline{\mathsf{V}}$.
(3.3.2) The self-enriched category $\underline{\mathsf{V}}$ is symmetric monoidal as a V-category.
(3.3.4) The standard enrichment of a symmetric monoidal functor is symmetric monoidal in the enriched sense.

(3.4.12) V**-Yoneda Bijection Theorem**. For each V-functor $F : \mathsf{C} \longrightarrow \underline{\mathsf{V}}$ and each $X \in \mathsf{C}$, there is a bijection of sets V-nat$(\mathcal{Y}^X, F) \cong \mathsf{V}(\mathbb{1}, FX)$. For each V-functor $G : \mathsf{C}^{\mathsf{op}} \longrightarrow \underline{\mathsf{V}}$ and each $Y \in \mathsf{C}$, there is a bijection of sets V-nat$(\mathcal{Y}_Y, G) \cong \mathsf{V}(\mathbb{1}, GY)$.
(3.5.1) A V-coend is initial among V-cowedges. A V-end is terminal among V-wedges.
(3.5.5) If V is cocomplete, then V-coends are computed by a coequalizer in V. If V is complete, then V-ends are computed by an equalizer in V.
(3.5.12) For V complete symmetric monoidal closed, the mapping object for V-functors to $\underline{\mathsf{V}}$ is given by a V-end.
(3.6.9) V**-Yoneda Lemma**. For a V-functor $F : \mathsf{C} \longrightarrow \mathsf{V}$ with C small, there is a V-natural isomorphism

$$F \xrightarrow{\cong} \mathsf{Map}(\mathcal{Y}^{(-)}, F).$$

(3.7.3) For V complete and cocomplete symmetric monoidal closed, the Day convolution and hom diagram are given by a V-coend and a V-end, respectively.
(3.7.8) V**-Yoneda Density Theorem**. For a V-functor $X : \mathcal{D} \longrightarrow \underline{\mathsf{V}}$ with $\mathcal{D}$ small, there is a V-natural isomorphism

$$\int^x \mathcal{D}(x, -) \otimes X_x \xrightarrow{\cong} X.$$

(3.7.13) There is an isomorphism $\mathsf{Map}(X, Y) \cong \mathsf{Hom}(X, Y)_e$, for V-diagrams $X$ and $Y$.
(3.7.22) **Day Convolution Theorem**. For a small symmetric monoidal V-category $\mathcal{D}$, the category of $\mathcal{D}$-shaped diagrams in V is symmetric monoidal closed with the Day convolution product, internal hom, and monoidal unit $J = \mathcal{Y}^e$.
(3.7.28) Precomposition with a symmetric monoidal V-functor induces a symmetric monoidal functor between diagram V-categories.
(3.8.4) Change of enrichment along a symmetric monoidal functor $U : \mathsf{V} \longrightarrow \mathsf{W}$ induces a symmetric monoidal functor from $\mathcal{D}$-V to $\mathcal{D}_U$-W.
(3.9.3) If C is tensored and cotensored over V, then $X \otimes -$ and $Y^{(-)}$ extend uniquely to V-functors that are V-adjoint to the respective co/represented V-functors.
(3.9.8) If $(F, U)$ is an adjunction of monoidal functors between symmetric monoidal closed categories and if $F^2$ is invertible, then $F$ transfers tensor and cotensor structure over its codomain to corresponding structure over its domain.
(3.9.15) The symmetric monoidal closed diagram category $\mathcal{D}$-V is enriched, tensored, and cotensored over V.

## **Chapter III.4**. Pointed Objects, Smash Products, and Pointed Homs

(4.1.6) Smash product with respect to a terminal object $T$ is given by a pushout from a monoidal product.
(4.2.1) Pointed hom with respect to a terminal object $T$ is given by a pullback from an internal hom.
(4.1.5, 4.1.8, 4.2.3) Suppose C is complete and cocomplete symmetric monoidal closed. Then $\mathsf{C}_*$ is complete and cocomplete symmetric monoidal closed with respect to the smash product and pointed hom.
(4.3.11) For a small symmetric monoidal category $\mathcal{D}$ with a null object, its pointed unitary enrichment over $(\mathsf{V}_*, \wedge, E)$ is given by taking wedge sums of $E$ over nonzero morphisms in $\mathcal{D}$.
(4.3.19) Assuming the basepoint of V is terminal and the basepoint of $\mathcal{D}$ is null, there is an equivalence of categories between pointed functors from $\mathcal{D}$ to $\mathsf{V}_*$ and $\mathsf{V}_*$-enriched functors from the pointed unitary enrichment of $\mathcal{D}$ to the self enrichment of $\mathsf{V}_*$.
(4.3.37) The category of pointed diagrams $\mathcal{D}_*$-$\mathsf{V}$ is complete and cocomplete symmetric monoidal closed. Moreover, it is enriched, tensored, and cotensored over $\mathsf{V}_*$.

**Chapter III.5**. Multicategories

(5.1.2) A multicategory has objects, $n$-ary operations, symmetric group actions, colored units, and composition that are subject to axioms for symmetry, associativity, unity, and equivariance.
(5.1.2) An operad is a multicategory with one object.
(5.1.11) Each small permutative category $\mathsf{C}$ has an endomorphism multicategory with the same objects and with $n$-ary operations given by morphisms out of $n$-fold sums in $\mathsf{C}$.
(5.1.12) A multifunctor satisfies axioms for symmetric group action, units, and composition.
(5.1.17) A multinatural transformation satisfies a naturality condition.
(5.1.20) There is a 2-category consisting of small multicategories, multifunctors, and multinatural transformations.
(5.1.21) The initial operad $\mathsf{I}$ has a single object and only one operation, which is the unit on its one object.
(5.2.1, 5.2.2) The terminal multicategory $\mathsf{T}$ has a single object and a single $n$-ary operation for each $n \geq 0$. The terminal multicategory is also known as the commutative operad, $\mathsf{Com}$.
(5.3.9) Taking endomorphism operads gives a 2-functor from $\mathsf{PermCat}^{\mathsf{su}}$ to $\mathsf{Multicat}_*$.
(5.4.1) A monad consists of an endofunctor together with multiplication and unit natural transformations such that the associativity and unit diagrams commute.
(5.4.2) A monad algebra consists of an object and a structure morphism such that associativity and unity diagrams commute.
(5.4.13) **Beck's Precise Tripleability Theorem**. An adjunction $L \dashv U$ is strictly monadic if and only if $U$ strictly creates coequalizers for parallel pairs $f, g$ for which $(Uf, Ug)$ has a split coequalizer.
(5.4.18) If $T$ is a monad on a complete and cocomplete category, and if $T$ preserves filtered colimits, then the category of $T$-algebras is complete and cocomplete.
(5.5.1) A multigraph consists of vertices and multiedges.
(5.5.4, 5.5.9) The forgetful functor from small multicategories to small multigraphs has a left adjoint.
(5.5.11) The category of small multicategories is strictly monadic over the category of small multigraphs.
(5.5.14) The category of small multicategories is complete and cocomplete.
(5.6.7) The sharp product of multicategories is generated by operations $\phi \times d$ and $c \times \psi$ subject to symmetry and compatibility axioms.
(5.6.12) The tensor product of multicategories is generated by those of the sharp product, and subject to an additional interchange relation.
(5.6.12) A multifunctor out of a tensor product of multicategories consists of an assignment on objects that is multifunctorial in each variable separately and that preserves the interchange relation.
(5.7.2, 5.7.4) The internal hom for multicategories has operations given by transformations that satisfy a naturality condition.

(5.7.14) The category of small multicategories is complete and cocomplete symmetric monoidal closed with monoidal product given by the tensor product and closed structure given by the internal hom.
(5.7.22) The category of small pointed multicategories is complete and cocomplete symmetric monoidal closed with monoidal product given by the smash product and closed structure given by the pointed hom.
(5.7.23) The symmetric monoidal structure on $\mathsf{Multicat}_*$ does not restrict along $\mathsf{End}$ to a symmetric monoidal structure on $\mathsf{PermCat}^{\mathsf{su}}$.

**Chapter III.6**. Enriched Multicategories

(6.1.1) A $\mathsf{V}$-enriched multicategory has $n$-ary operation objects, symmetric group action, colored units, and composition morphisms in $\mathsf{V}$. These satisfy axioms given by commutative diagrams in $\mathsf{V}$ for symmetric group action, associativity, unity, and equivariance.
(6.1.8) A $\mathsf{V}$-enriched operad is a $\mathsf{V}$-multicategory with one object.
(6.1.9) An object of an enriched multicategory has a $\mathsf{V}$-enriched endomorphism operad.
(6.1.10) A $\mathsf{V}$-enriched multifunctor satisfies axioms given by commutative diagrams in $\mathsf{V}$ for symmetric group action, units, and composition.
(6.1.14) An algebra $c$ over a $\mathsf{V}$-enriched operad $\mathsf{P}$ is given by a $\mathsf{V}$-enriched operad morphism $\mathsf{P} \longrightarrow \mathsf{End}(c)$.
(6.1.15) A $\mathsf{V}$-enriched multinatural transformation satisfies a $\mathsf{V}$-naturality diagram in $\mathsf{V}$.
(6.1.18) There is a 2-category consisting of small $\mathsf{V}$-enriched multicategories together with $\mathsf{V}$-enriched multifunctors and multinatural transformations.
(6.2.9) For a symmetric monoidal functor $U : \mathsf{V} \longrightarrow \mathsf{W}$, change of enrichment along $U$ provides a 2-functor from small $\mathsf{V}$-enriched multicategories to small $\mathsf{W}$-enriched multicategories.
(6.3.3, 6.3.6) A symmetric monoidal $\mathsf{V}$-category has a $\mathsf{V}$-enriched endomorphism multicategory with $\mathsf{V}$-objects of $n$-ary operations given by morphism objects out of $n$-fold left normalized products in $\mathsf{K}$.
(6.3.10) A symmetric monoidal $\mathsf{V}$-functor induces a $\mathsf{V}$-enriched multifunctor between $\mathsf{V}$-enriched endomorphism multicategories of its domain and codomain.
(6.4.3) The tensor product is a $\mathsf{Cat}$-enriched symmetric monoidal product for the 2-category of small multicategories.
(6.4.4) The smash product is a $\mathsf{Cat}$-enriched symmetric monoidal product for the 2-category of small pointed multicategories.
(6.4.5) Each of $\mathsf{Multicat}$ and $\mathsf{Multicat}_*$ has the structure of a $\mathsf{Cat}$-enriched multicategory induced by the tensor and smash product, respectively.
(6.5.1) The $\mathsf{Cat}$-enriched multicategory structure on $\mathsf{Multicat}_*$ induces a corresponding structure on $\mathsf{PermCat}^{\mathsf{su}}$.
(6.5.4) Multilinear functors of permutative categories consist of functors out of a Cartesian product together with linearity constraints. They are subject to axioms for unity, constraint unity, constraint associativity, constraint symmetry, and constraint 2-by-2.
(6.5.11) Multilinear transformations between multilinear functors satisfy multilinearity conditions with respect to linearity constraints and identities.
(6.5.10, 6.5.13) The categories of $n$-ary operations in $\mathsf{PermCat}^{\mathsf{su}}$ are canonically isomorphic to the corresponding categories of $n$-linear functors and transformations.

(6.6.13) The Cat-enriched multicategory structure of $\mathsf{PermCat}^{\mathsf{su}}$ is described explicitly in terms of multilinear functors and transformations.

## **Part** III.2. Algebraic $K$-Theory

### **Chapter** III.7. Homotopy Theory Background

(7.1.16) The geometric realization of the standard $n$-simplex is the topological $n$-simplex.
(7.1.19) The category of simplicial sets is symmetric monoidal closed with the monoidal product given by the levelwise Cartesian product.
(7.1.23) The category of pointed simplicial sets is symmetric monoidal closed with monoidal product given by the levelwise smash product.
(7.2.4) The nerve of a small category is a simplicial set with $p$-simplices given by strings of $p$ composable arrows.
(7.2.5 (1)) A natural transformation between functors induces a simplicial homotopy on nerves.
(7.2.5 (2)) An adjunction of functors induces a simplicial homotopy equivalence on nerves.
(7.3.7) The category of symmetric sequences is symmetric monoidal closed with monoidal product given by Day convolution.
(7.4.5) The category of symmetric spectra is the category of left modules over the symmetric sphere.
(7.4.6) A symmetric spectrum consists of a symmetric sequence with structure maps satisfying unity, associativity, and equivariance axioms.
(7.5.5) The category of symmetric spectra is complete and cocomplete.
(7.6.1) The smash product of symmetric spectra is given by a coequalizer of actions by $S$.
(7.6.8) The internal hom for symmetric spectra is given by an equalizer of actions by $S$.
(7.6.15) The category of symmetric spectra is complete and cocomplete symmetric monoidal closed.
(7.8.8) Every level equivalence of symmetric spectra is a stable equivalence.

### **Chapter** III.8. Segal $K$-Theory of Permutative Categories

(8.1.8) A $\Gamma$-object in a pointed category $\mathsf{C}$ is a pointed functor from $\mathcal{F}$ to $\mathsf{C}$.
(8.2.6) The construction $\mathsf{K}^{\mathcal{F}}$ is a functor from $\Gamma$-simplicial sets to symmetric spectra.
(8.3.13) For a small permutative category $\mathsf{C}$, there are three variant constructions of $\Gamma$-categories, $\mathsf{C}^{\mathcal{F}} = \mathsf{C}^{\mathcal{F}}_{\cong}$, $\mathsf{C}^{\mathcal{F}}_{\mathrm{lax}}$, and $\mathsf{C}^{\mathcal{F}}_{\mathrm{co}}$.
(8.3.21) For a small permutative category $\mathsf{C}$, each of the $\Gamma$-simplicial sets $N_*\mathsf{C}^{\mathcal{F}}_{\cong}$, $N_*\mathsf{C}^{\mathcal{F}}_{\mathrm{lax}}$, and $N_*\mathsf{C}^{\mathcal{F}}_{\mathrm{co}}$ is special, and all three are levelwise weakly-equivalent.
(8.4.5) The partition multicategory $\mathcal{M}\underline{1}$ has two objects, $\varnothing$ and $\{1\}$, with operations given by partitions.
(8.4.7) The partition multicategory $\mathcal{M}$ defines a pointed functor from $\mathcal{F}^{\mathsf{op}}$ to $\mathsf{Multicat}_*$.
(8.4.8, 8.4.10) For a small permutative category $\mathsf{C}$, there is an isomorphism of $\Gamma$-categories $\mathsf{J}^{\mathsf{Se}}\mathsf{C} \cong \mathsf{C}^{\mathcal{F}}_{\mathrm{lax}}$.

### **Chapter** III.9. Categories of $\mathcal{G}_*$-Objects

(9.1.7) The objects of $\mathcal{G}$ are tuples of objects of $\mathcal{F}$ subject to certain basepoint identifications.

(9.1.15) Smash product of pointed finite sets defines a strict symmetric monoidal functor from $\mathcal{G}$ to $\mathcal{F}$.
(9.2.1) A $\mathcal{G}_*$-object in a pointed category $\mathsf{C}$ is a pointed functor from $\mathcal{G}$ to $\mathsf{C}$.
(9.2.15) If $\mathsf{C}$ is complete and cocomplete symmetric monoidal closed with terminal basepoint, then the category for $\mathcal{G}_*$-objects in $\mathsf{C}$ is complete and cocomplete symmetric monoidal closed with monoidal product given by Day convolution.
(9.2.19) The nerve induces a symmetric monoidal $\mathsf{sSet}_*$-functor from small $\mathcal{G}_*$-categories to $\mathcal{G}_*$-simplicial sets.
(9.3.16) The construction $\mathsf{K}^{\mathcal{G}}$ is a functor from $\mathcal{G}_*$-simplicial sets to symmetric spectra whose restriction along $\wedge^*$ is equal to $\mathsf{K}^{\mathcal{F}}$.
(9.4.9) The functor $\mathsf{K}^{\mathcal{G}}$ is a unital symmetric monoidal $\mathsf{sSet}$-functor.
(9.4.18) The functors $\wedge^*$ and $\mathsf{K}^{\mathcal{F}}$ are symmetric monoidal $\mathsf{sSet}_*$-functors.

## Chapter III.10. Elmendorf-Mandell $K$-Theory of Permutative Categories

(10.1.6) The partition multicategory $\mathcal{M}$ is a symmetric monoidal functor, with monoidal constraint given by the partition product $\prod$.
(10.1.12) The partition products $\prod_{\underline{1},b}$ and $\prod_{b,\underline{1}}$ are isomorphisms.
(10.1.25) The category of left $\mathcal{M}\underline{1}$-modules in $\mathsf{Multicat}_*$ is a full subcategory and the smash product over $\mathcal{M}\underline{1}$ is isomorphic to that of $\mathsf{Multicat}_*$.
(10.1.33) The category of left $\mathcal{M}\underline{1}$-modules is symmetric monoidal in the $\mathsf{Cat}_*$-enriched sense.
(10.1.35) The category of left $\mathcal{M}\underline{1}$-modules is complete and cocomplete.
(10.2.7) If $\mathsf{C}$ is a small permutative category, $\mathsf{End}(\mathsf{C})$ has a canonical left $\mathcal{M}\underline{1}$-module structure. Taking this structure, $\mathsf{End}$ factors through $\mathsf{Mod}^{\mathcal{M}\underline{1}}$.
(10.2.8, 10.2.9) The category of left $\mathcal{M}\underline{1}$-modules is a proper subcategory of $\mathsf{Multicat}_*$. In particular, the monoidal unit $\mathsf{S}$ is not an $\mathcal{M}\underline{1}$-module.
(10.2.10) The symmetric monoidal structure on $\mathsf{Mod}^{\mathcal{M}\underline{1}}$ does not restrict to $\mathsf{PermCat}^{\mathsf{su}}$.
(10.3.3) The smash product of partition multicategories, $\mathcal{T}$, defines a pointed functor from $\mathcal{G}^{\mathsf{op}}$ to $\mathsf{Mod}^{\mathcal{M}\underline{1}}$.
(10.3.7) The functor $\mathcal{T}$ is strictly unital strong symmetric monoidal.
(10.3.13) The monoidal constraint for the partition $J$-theory $\mathsf{J}^{\mathcal{T}}$ uses the inverse monoidal constraint for $\mathcal{T}$.
(10.3.17) The partition $J$-theory $\mathsf{J}^{\mathcal{T}}$ is a symmetric monoidal $\mathsf{Cat}_*$-functor.
(10.3.25) Elmendorf-Mandell $J$-theory $\mathsf{J}^{\mathsf{EM}} = \mathsf{J}^{\mathcal{T}} \circ \mathsf{End}$ is a $\mathsf{Cat}$-enriched multifunctor.
(10.3.32) Elmendorf-Mandell $K$-theory $\mathsf{K}^{\mathsf{EM}} = \mathsf{K}^{\mathcal{G}} N_* \mathsf{J}^{\mathsf{EM}}$ is a $\mathsf{sSet}$-enriched multifunctor.
(10.3.33) The multifunctor $\mathsf{K}^{\mathsf{EM}}$ preserves enriched operad actions.
(10.4.18) For a small permutative category $\mathsf{C}$ there are three variant constructions of $\mathcal{G}_*$-categories, $\mathsf{C}^{\mathcal{G}} = \mathsf{C}^{\mathcal{G}}_{\mathrm{lax}}$, $\mathsf{C}^{\mathcal{G}}_{\cong}$, and $\mathsf{C}^{\mathcal{G}}_{\mathrm{co}}$.
(10.5.1) For a small permutative category $\mathsf{C}$, there is an isomorphism of $\mathcal{G}_*$-categories $\mathsf{J}^{\mathsf{EM}}\mathsf{C} \cong \mathsf{C}^{\mathcal{G}}_{\mathrm{lax}}$.
(10.6.10) There is a level equivalence of symmetric spectra $\mathsf{K}^{\mathsf{Se}}\mathsf{C} \longrightarrow \mathsf{K}^{\mathsf{EM}}\mathsf{C}$ for each small permutative category $\mathsf{C}$. It is natural with respect to strictly unital symmetric monoidal functors.
(10.7.16) There is a level equivalence $\mathsf{K}^{\mathsf{EM}}_{\cong}\mathsf{C} \longrightarrow \mathsf{K}^{\mathsf{EM}}_{\mathsf{lax}}$ for each small permutative category $\mathsf{C}$.

(10.7.19 (4)) There is a level equivalence $\mathsf{K}^{\mathsf{EM}}_{\cong} \longrightarrow \mathsf{K}^{\mathsf{EM}}_{\mathsf{co}}\mathsf{C}$ for each small permutative category $\mathsf{C}$.
(10.7.22, 10.7.27) The $\mathcal{G}_*$-category morphisms $\mathsf{C}^{\mathcal{G}}_{\cong} \longrightarrow \mathsf{C}^{\mathcal{G}}$ and $\mathsf{C}^{\mathcal{G}}_{\cong} \longrightarrow \mathsf{C}^{\mathcal{G}}_{\mathsf{co}}$ are components of $\mathsf{Cat}$-enriched multinatural transformations.

**Chapter III.11**. *K*-Theory of Ring and Bipermutative Categories

(11.1.4) $\mathsf{As}$ is an operad.
(11.1.7) $\mathsf{As}$ is generated by $\mathrm{id}_0$ and $\mathrm{id}_2$, which are subject to unity and associativity relations.
(11.1.15) $\mathsf{As}$ is the operad for monoids.
(11.2.16) $\mathsf{As}$ detects ring category structures on small permutative categories.
(11.3.2) A strict ring symmetric spectrum is a symmetric spectrum equipped with multiplication and unit morphisms in symmetric sequences that satisfy compatibility, associativity, and unity axioms.
(11.3.13–11.3.15) The sphere spectrum, the suspension spectrum of a monoid in pointed simplicial sets, and the Eilenberg-Mac Lane spectrum of a ring are strict ring symmetric spectra.
(11.3.16) $\mathsf{K}^{\mathsf{EM}}\mathsf{C}$ is a strict ring symmetric spectrum for each small ring category $\mathsf{C}$.
(11.3.17) For each small permutative category, the *K*-theory of its (tight) endomorphism ring category is a strict ring symmetric spectrum.
(11.3.18) The *K*-theory of the additive distortion category is a strict ring symmetric spectrum.
(11.3.19) For each small tight bimonoidal category, the *K*-theory of its right/left rigid bimonoidal strictification is a strict ring symmetric spectrum.
(11.4.7) The translation category functor $E$ is a right adjoint.
(11.4.11) Each morphism in the Barratt-Eccles operad decomposes into a categorical composite of $\phi v$ with $v$ a permutation and $\phi$ an operadic composite of one $\tau : \mathrm{id}_2 \longrightarrow (1,2)$ and identity morphisms.
(11.4.14) The Barratt-Eccles operad is generated by two objects and one isomorphism, which are subject to relations that are formally identical to those of a permutative category.
(11.4.26) The Barratt-Eccles operad is the $\mathsf{Cat}$-enriched operad for permutative categories.
(11.5.5) The Barratt-Eccles operad detects bipermutative category structures on small permutative categories.
(11.6.3) The Barratt-Eccles operad is an $E_\infty$-operad.
(11.6.6 (1)) Each commutative monoid in $\mathsf{SymSp}$ is an $E_\infty$-symmetric spectrum.
(11.6.6 (2)) Each $E_\infty$-symmetric spectrum via the Barratt-Eccles operad has a strict ring structure.
(11.6.7) The symmetric sphere is a commutative monoid in $\mathsf{SymSp}$.
(11.6.9) The suspension spectrum of a commutative monoid in $\mathsf{sSet}_*$ is a commutative monoid in $\mathsf{SymSp}$.
(11.6.10) The Eilenberg-Mac Lane spectrum of a commutative ring is a commutative monoid in $\mathsf{SymSp}$.
(11.6.12) $\mathsf{K}^{\mathsf{EM}}\mathsf{C}$ is an $E_\infty$-symmetric spectrum for each small bipermutative category $\mathsf{C}$.
(11.6.13) The *K*-theory of each small right/left bipermutative category is an $E_\infty$-symmetric spectrum. For example, this applies to the finite ordinal category $\Sigma$, its variant $\Sigma'$, $\mathsf{Vect}^{\mathbb{C}}_{\mathsf{c}}$, and the distortion category.

(11.6.14) For each small tight symmetric bimonoidal category, the $K$-theory of its right/left bipermutative strictification is an $E_\infty$-symmetric spectrum. For example, this applies to small distributive symmetric monoidal categories, the symmetric bimonoidal groupoid $\Pi$, and the bimonoidal symmetric center of a small tight braided bimonoidal category.
(11.6.15) For each small braided ring category whose left factorization morphism is a natural epimorphism, the $K$-theory of its symmetric center is an $E_\infty$-symmetric spectrum.

## **Chapter III.12**. *K*-Theory of Braided Ring Categories

(12.1.10) The braid operad is a $\mathsf{Cat}$-enriched operad.
(12.2.4) The braid operad is an $E_2$-operad.
(12.3.6) Each morphism in the braid operad admits a categorical decomposition into isomorphisms of the form $\phi v$ with $v$ a permutation and $\phi$ an operadic composite of one $s_1^{\pm 1} : \mathrm{id}_2 \longrightarrow (1,2)$ and identity morphisms.
(12.3.10) The braid operad is generated by two objects and one isomorphism, which are subject to relations that are formally identical to those of a braided strict monoidal category.
(12.3.22) The braid operad is the $\mathsf{Cat}$-enriched operad for braided strict monoidal categories.
(12.4.5) The braid operad detects braided ring category structures on small permutative categories.
(12.5.2 (1)) An $E_\infty$-symmetric spectrum via the Barratt-Eccles operad has an $E_2$-structure.
(12.5.2 (2)) An $E_2$-symmetric spectrum via the braid operad has a strict ring structure.
(12.5.3) $\mathsf{K}^{\mathsf{EM}}\mathsf{C}$ is an $E_2$-symmetric spectrum for each small braided ring category $\mathsf{C}$.
(12.5.4) The $K$-theory of the braided distortion category is an $E_2$-symmetric spectrum.
(12.5.5) For each small tight ring category, the $K$-theory of its bimonoidal Drinfeld center is an $E_2$-symmetric spectrum.
(12.5.6) For each small tight braided bimonoidal category, the $K$-theory of its right/left permbraided strictification is an $E_2$-symmetric spectrum. For example, this applies to a small abelian category with a compatible braided monoidal structure, Fibonacci anyons, Ising anyons, and the bimonoidal Drinfeld center of a small tight bimonoidal category.

## **Chapter III.13**. *K*-Theory of $E_n$-Monoidal Categories

(13.1.20) $\mathsf{Mon}^n$ is a $\mathsf{Cat}$-enriched operad.
(13.1.23) $\mathsf{Mon}^1$ is the associative operad.
(13.2.1) $\mathsf{Mon}^n$ is an $E_n$-operad.
(13.3.3) $\mathsf{Mon}^n$ is generated by the objects $\mathbb{1}$ and $\{1 \otimes_i 2\}_{i=1}^n$ and the exchange morphisms $\{\eta_{1,2,3,4}^{i,j}\}_{1\le i<j\le n}$, which are subject to relations that are formally identical to those of an $n$-fold monoidal category.
(13.3.18) $\mathsf{Mon}^n$ is the $\mathsf{Cat}$-enriched operad for $n$-fold monoidal categories.
(13.4.12) $\mathsf{Mon}^n$ detects $E_n$-monoidal category structures on small permutative categories.

(13.5.1) The canonical $\mathsf{Cat}$-enriched operad morphism $\mathsf{As} \longrightarrow E\mathsf{As}$ factors through $\mathsf{Mon}^n$. As a result, an $E_\infty$-symmetric spectrum via the Barratt-Eccles operad induces an $E_n$-structure. An $E_{n+1}$-structure via $\mathsf{Mon}^{n+1}$ induces an $E_n$-structure.
(13.5.2) $\mathsf{K}^{\mathsf{EM}}\mathsf{C}$ is an $E_n$-symmetric spectrum for each small $E_n$-monoidal category $\mathsf{C}$.
(13.5.3) For each small category, the $K$-theory of its free $E_n$-monoidal category is an $E_n$-symmetric spectrum.

# List of Notations

## Part I.1

| Chapter I.1 | Page | Description |
|---|---|---|
| $\mathsf{Ob}(\mathsf{C})$, $\mathsf{Ob}\mathsf{C}$ | I.7 | objects in a category $\mathsf{C}$ |
| $\mathsf{C}(X,Y)$, $\mathsf{C}(X;Y)$ | I.8 | set of morphisms $X \longrightarrow Y$ |
| $1_X$ | I.8 | identity morphism |
| $g \circ f$, $gf$ | I.8 | composition of morphisms |
| $\cong$, $\xrightarrow{\cong}$ | I.8 | an isomorphism |
| $F : \mathsf{C} \longrightarrow \mathsf{D}$ | I.8 | a functor |
| $\mathrm{Id}_{\mathsf{C}}$, $1_{\mathsf{C}}$ | I.9 | identity functor |
| $\mathbf{1}$ | I.9 | terminal category |
| $\theta_X$ | I.9 | a component of a natural transformation $\theta$ |
| $1_F$ | I.9 | identity natural transformation |
| $\phi\theta$ | I.9 | vertical composition |
| $\theta' * \theta$ | I.9 | horizontal composition |
| $(L, R, \phi)$, $L \dashv R$ | I.10 | an adjunction |
| $\eta$, $\varepsilon$ | I.10 | unit and counit of an adjunction |
| $\mathsf{colim}\, F$ | I.11 | colimit |
| $\mathsf{lim}\, F$ | I.11 | limit |
| $\varnothing$, $\varnothing^{\mathsf{C}}$ | I.12 | an initial object |
| $\coprod$, $\amalg$ | I.12 | a coproduct |
| $\otimes$ | I.14 | monoidal product |
| $\mathbb{1}$ | I.14 | monoidal unit |
| $\alpha$ | I.14 | associativity isomorphism |
| $\lambda$, $\rho$ | I.15 | unit isomorphisms |
| $(X, \mu, \eta)$ | I.15 | a monoid |
| $(Y, \Delta, \varepsilon)$ | I.16 | a comonoid |
| $(F, F^2, F^0)$ | I.16 | a monoidal functor |
| $\xi$ | I.18 | symmetry isomorphism |
| $(\mathsf{Set}, \times, *)$ | I.19 | category of sets |
| $(\mathsf{Cat}, \times, \mathbf{1})$ | I.19 | category of small categories |
| $(\mathsf{Vect}^{\Bbbk}, \otimes, \Bbbk)$ | I.19 | category of $\Bbbk$-vector spaces |
| $[-,-]$ | I.19 | internal hom |
| $e$, $-$, $u \,\square\, v$ | I.19 | words |
| $\sigma(-)$ | I.20 | a left permutation |
| $w\sigma$ | I.20 | a permuted word |

| Chapter I.2 | | |
|---|---|---|
| $(\oplus, 0, \alpha^{\oplus}, \lambda^{\oplus}, \rho^{\oplus}, \xi^{\oplus})$ | I.25 | additive structure |
| $(\otimes, \mathbb{1}, \alpha^{\otimes}, \lambda^{\otimes}, \rho^{\otimes}, \xi^{\otimes})$ | I.25 | multiplicative structure |
| $\lambda^{\bullet}$, $\rho^{\bullet}$ | I.25 | multiplicative zeros |

## Chapter I.7

## Chapter I.8

| | | |
|---|---|---|
| $(a^{\boxtimes}, a^{\boxtimes\bullet}, \eta^a, \varepsilon^a)$ | I.383 | monoidal associator |
| $(\ell^{\boxtimes}, \ell^{\boxtimes\bullet}, \eta^\ell, \varepsilon^\ell)$ | I.387 | left monoidal unitor |
| $(r^{\boxtimes}, r^{\boxtimes\bullet}, \eta^r, \varepsilon^r)$ | I.390 | right monoidal unitor |
| $\pi$ | I.392 | pentagonator |
| $\mu$ | I.400 | middle 2-unitor |
| $\lambda^{\boxtimes}$ | I.402 | left 2-unitor |
| $\rho^{\boxtimes}$ | I.404 | right 2-unitor |
| ${}_\theta A$ | I.409 | row permutation of $A$ by $\theta$ |
| $A^\sigma$ | I.409 | column permutation of $A$ by $\sigma$ |
| $\mathbb{1}^\sigma$ | I.409 | permutation matrix of $\sigma$ |
| $r^\sigma_A$ | I.410 | natural isomorphism $A\mathbb{1}^\sigma \xrightarrow{\cong} A^\sigma$ |
| $\ell^\theta_A$ | I.411 | natural isomorphism $\mathbb{1}^\theta A \xrightarrow{\cong} {}_{\theta^{-1}}A$ |
| $\tau_{m,n}$ | I.412 | permutation in $\Sigma_{mn}$ that transposes an $n \times m$ matrix |
| $(\beta, \beta^\bullet, \eta^\beta, \varepsilon^\beta)$ | I.418 | braiding |
| $h_{m\|n,p}$ | I.420 | comparison 2-cell for $R_{-\|--}$ |
| $R_{-\|--}$ | I.421 | left hexagonator |
| $h_{m,n\|p}$ | I.424 | comparison 2-cell for $R_{--\|-}$ |
| $R_{--\|-}$ | I.425 | right hexagonator |
| $\nu$ | I.428 | syllepsis |

## Chapter I.9

| | | |
|---|---|---|
| $\odot$ | I.437 | Gray monoid multiplication on $\mathsf{Mat}^{\mathsf{C}}$ |
| $\Sigma_{A,B}$ | I.438 | transition 2-cells |
| $\beta$ | I.446 | Gray symmetry |
| $\boxdot$ | I.450 | permutative 2-category multiplication on $\mathsf{Mat}^{\mathsf{C}}$ |

## Part II.1

## Chapter II.1

| | | |
|---|---|---|
| $B_n$ | II.8 | braid group on $n$ strings |
| $s_1, \ldots, s_{n-1}$ | II.8 | generating braids in $B_n$ |
| $s_i^{(n)}$ | II.8 | $s_i$ in $B_n$ |
| id, $\mathrm{id}_n$ | II.8 | identity braid in $B_n$ |
| $\mathcal{I}$ | II.8 | unit interval $[0,1]$ |
| $\oplus$ | II.9 | sum braid |
| $\pi(b), \overline{b}$ | II.10 | underlying permutation of $b$ |
| $\sigma\langle \underline{k} \rangle$ | II.10 | block permutation induced by $\sigma \in \Sigma_n$ |
| $b\langle \underline{k} \rangle$ | II.11 | block braid induced by $b \in B_n$ |
| $\tau\langle m,n \rangle$ | II.13 | interval-swapping permutation in $\Sigma_{m+n}$ |
| $b^{\oplus}_{m,n}$ | II.14 | elementary block braid induced by $s_1 \in B_2$ |
| $(\overline{\mathsf{C}}, \overline{\otimes}, \overline{\mathbb{1}}, \overline{\alpha}, \overline{\lambda}, \overline{\rho}, \overline{\xi})$ | II.26 | Drinfeld center of $\mathsf{C}$ |
| $(A; \beta^A)$ | II.26 | an object in $\overline{\mathsf{C}}$ |
| $\mathsf{C}^{\mathsf{sym}}$ | II.35 | symmetric center of $\mathsf{C}$ |
| $\mathsf{br}(\phi)$ | II.37 | underlying braid of $\phi$ |

## Chapter II.2

| | | |
|---|---|---|
| $(\oplus, 0, \alpha^\oplus, \lambda^\oplus, \rho^\oplus, \xi^\oplus)$ | II.45 | additive structure |
| $(\otimes, \mathbb{1}, \alpha^\otimes, \lambda^\otimes, \rho^\otimes, \xi^\otimes)$ | II.45 | multiplicative structure |
| $\lambda^\bullet, \rho^\bullet$ | II.46 | multiplicative zeros |
| $\delta^l, \delta^r$ | II.46 | distributivity morphisms |
| $\mathsf{Ab}$ | II.50 | category of abelian groups |
| $0 : A \longrightarrow B$ | II.51 | zero morphism |
| $i_1, i_2$ | II.51 | inclusions |
| $p_1$, $p_2$ | II.51 | projections |
| $A \oplus B$ | II.51 | direct sum of objects $A$ and $B$ |
| $f \oplus f'$ | II.52 | direct sum morphism |

## Chapter II.3

| | | |
|---|---|---|
| $\mathsf{Vect}^{\Bbbk}$ | II.71 | category of $\Bbbk$-vector spaces |
| $A^{\otimes n}$ | II.71 | $A \otimes \cdots \otimes A$ with $n$ copies of $A$ and $A^{\otimes 0} = \Bbbk$ |
| $A^{\oplus n}$ | II.71 | $A \oplus \cdots \oplus A$ with $n$ copies of $A$ and $A^{\oplus 0} = 0$ |
| $\mathsf{Vect}^{\Bbbk}_{\otimes}$ | II.71 | $\mathsf{Vect}^{\Bbbk}$ with the tensor product |
| $\mathsf{Vect}^{\Bbbk}_{\oplus}$ | II.71 | $\mathsf{Vect}^{\Bbbk}$ with the direct sum |
| $(A, \mu, \eta, \Delta, \varepsilon)$ | II.72 | a bialgebra |
| $\sum_i s_i' \otimes s_i''$ | II.73 | an element in $A^{\otimes 2}$ |
| $S^{\mathsf{op}}$ | II.73 | $\xi^{\otimes} S$ |
| $S_{12}$, $S_{13}$, $S_{23}$ | II.73 | elements obtained from $S \in A^{\otimes 2}$ by inserting 1 |
| $\sum_{(x)} x^{(1)} \otimes x^{(2)}$ | II.74 | Sweedler's notation for comultiplication |
| $\Delta^{\mathsf{op}}$ | II.74 | opposite comultiplication $\xi^{\otimes}\Delta$ |
| $\Bbbk G$ | II.75 | group bialgebra of $G$ |
| $U\mathfrak{g}$ | II.76 | universal enveloping bialgebra of $\mathfrak{g}$ |
| $H_4$ | II.77 | Sweedler's 4-dimensional bialgebra |
| $\mathbb{C}\mathbb{Z}_n$ | II.78 | anyonic quantum groups |
| $\mathsf{Mod}(A)_{\otimes}$ | II.80 | $\mathsf{Mod}(A)$ with the tensor product |
| $\mathsf{Mod}(A)_{\oplus}$ | II.83 | $\mathsf{Mod}(A)$ with the direct sum |
| $\mathsf{Vect}^{\mathbb{C}}_{\mathsf{sk}}$ | II.85 | abelian category with linear maps $\mathbb{C}^m \longrightarrow \mathbb{C}^n$ as morphisms |
| $\mathcal{F}^{\mathsf{any}}$ | II.86 | Fibonacci anyons |
| $0 = (0;0)$ | II.86 | additive zero in $\mathcal{F}^{\mathsf{any}}$ |
| $\mathbb{1} = (1;0)$ | II.86 | vacuum in $\mathcal{F}^{\mathsf{any}}$ |
| $\tau = (0;1)$ | II.86 | non-abelian anyon in $\mathcal{F}^{\mathsf{any}}$ |
| $q$ | II.90 | reciprocal of the golden ratio |
| $z$ | II.92 | $e^{3\pi i/5}$ |
| $\mathcal{I}^{\mathsf{any}}$ | II.95 | Ising anyons |
| $0 = (0;0;0)$ | II.95 | additive zero in $\mathcal{I}^{\mathsf{any}}$ |
| $\mathbb{1} = (1;0;0)$ | II.95 | vacuum in $\mathcal{I}^{\mathsf{any}}$ |
| $\sigma = (0;1;0)$ | II.95 | non-abelian anyon in $\mathcal{I}^{\mathsf{any}}$ |
| $\psi = (0;0;1)$ | II.95 | fermion in $\mathcal{I}^{\mathsf{any}}$ |
| $w$ | II.105 | $e^{\pi i/8}$ |

## Chapter II.4

| | | |
|---|---|---|
| $\overline{\mathsf{C}}^{\mathsf{bi}}$ | II.113 | bimonoidal Drinfeld center |
| $(\overline{\oplus}, \overline{0}, \alpha^{\overline{\oplus}}, \lambda^{\overline{\oplus}}, \rho^{\overline{\oplus}}, \xi^{\overline{\oplus}})$ | II.113 | additive structure of $\overline{\mathsf{C}}^{\mathsf{bi}}$ |
| $(\overline{\otimes}, \overline{\mathbb{1}}, \alpha^{\overline{\otimes}}, \lambda^{\overline{\otimes}}, \rho^{\overline{\otimes}}, \xi^{\overline{\otimes}})$ | II.113 | multiplicative structure of $\overline{\mathsf{C}}^{\mathsf{bi}}$ |
| $\overline{\lambda}^{\bullet}, \overline{\rho}^{\bullet}$ | II.113 | multiplicative zeros of $\overline{\mathsf{C}}^{\mathsf{bi}}$ |
| $\overline{\delta}^{l}, \overline{\delta}^{r}$ | II.113 | distributivity morphisms of $\overline{\mathsf{C}}^{\mathsf{bi}}$ |
| $(A;\beta^A)$ | II.114 | an object in $\overline{\mathsf{C}}^{\mathsf{bi}}$ |
| $\mathsf{C}^{\mathsf{sym}}$ | II.127 | bimonoidal symmetric center |

## Chapter II.5

| | | |
|---|---|---|
| $\mathcal{D}^{\mathsf{br}}$ | II.137 | braided distortion category |
| $\underline{r}$ | II.137 | an object in $\mathcal{D}^{\mathsf{br}}$ |
| $\varnothing$ | II.137 | empty sequence in $\mathcal{D}^{\mathsf{br}}$ |
| $\underline{\sigma}$ | II.137 | a morphism in $\mathcal{D}^{\mathsf{br}}$ |
| $\int_{\mathsf{C}} F$ | II.156 | Grothendieck construction of $F : \mathsf{C}^{\mathsf{op}} \longrightarrow \mathsf{Cat}$ |
| $\mathcal{B}$ | II.157 | braid category |

## Chapter II.6

## Chapter II.7

## Chapter II.8

## Part II.2
## Chapter II.9

## Chapter III.2

## Chapter III.3

## Chapter III.4

## Chapter III.5

## Chapter III.6

## Part III.2

## Chapter III.7

| | | |
|---|---|---|
| $BC = \lvert NC \rvert$ | III.272 | classifying space of small category $\mathsf{C}$ |
| $\Sigma$ | III.272 | finite ordinal category |
| $X : \Sigma \longrightarrow \mathsf{sSet}_*$ | III.272 | symmetric sequence |
| $X \square Y$ | III.273 | Day convolution of symmetric sequences |
| $I, \lambda, \rho, \alpha, \xi$ | III.273 | monoidal data for $\square$ |
| $\Sigma[n] = \Sigma(\underline{n}, -)$ | III.273 | corepresented symmetric sequence |
| $X \wedge A$, $X^A$ | III.274 | simplicial tensor and cotensor for symmetric sequence $X$ |
| $\mathsf{Map}_\Sigma(X,Y)$, $\mathsf{Hom}_\Sigma(X,Y)$ | III.274 | symmetric mapping object and symmetric hom object |
| $i_n \dashv \mathsf{ev}_n$ | III.274 | evaluation at $\underline{n}$ for symmetric sequences, and left adjoint $i_n$ |
| $X^{\wedge k}$ | III.275 | right normalized iterated smash product |
| $\mathsf{Sym}(X)$ | III.275 | symmetric sequence given by smash powers of $X \in \mathsf{sSet}_*$ |
| $S = \mathsf{Sym}(S^1)$ | III.276 | symmetric sphere |
| $(A, \mu, \eta)$, $(M, \theta)$ | III.276 | commutative monoid and left module |
| $\mathsf{SymSp}$ | III.276 | category of symmetric spectra |
| $\rho, \rho^{\Sigma}_{p,q}, \rho_{p,q}$ | III.276 | structure morphisms for symmetric spectrum |
| $\Sigma^\infty K$ | III.277 | suspension spectrum of $K \in \mathsf{sSet}_*$ |
| $HR$ | III.278 | Eilenberg-Mac Lane spectrum of commutative ring $R$ |
| $(A \otimes -, \pi, \nu)$ | III.279 | monad associated to monoid $A$ |
| $X \wedge Y = X \square_S Y$ | III.282 | smash product of symmetric spectra |
| $\mathsf{Hom}_S(X,Y)$ | III.284 | internal hom for symmetric spectra |
| $F_n \dashv \mathsf{ev}_n$ | III.285 | evaluation at $\underline{n}$ for symmetric spectra, and left adjoint $F_n$ |
| $\mathsf{Map}_S(X,Y)$ | III.287 | mapping simplicial set for symmetric spectra $X$ and $Y$ |
| $\mathsf{End}(X)$ | III.287 | endomorphism simplicial operad for symmetric spectrum $X$ |
| $(\mathcal{L}, \mathcal{R})$ | III.287 | weak factorization system |
| $\mathsf{M}, \mathcal{W}, \mathcal{C}, \mathcal{F}$ | III.288 | model category, weak equivalences, cofibrations, and fibrations |
| $A^c$, $A^f$ | III.288 | cofibrant replacement and fibrant replacement |
| $\mathsf{Cyl}(A)$, $\mathsf{Path}(A)$ | III.288 | cylinder object and path object |
| $\gamma : \mathsf{M} \longrightarrow Ho\mathsf{M}$ | III.290 | homotopy category and canonical functor of $\mathsf{M}$ |
| $\mathsf{M} \simeq_Q \mathsf{N}$ | III.292 | Quillen equivalence between model categories |
| $\mathsf{Cell}(\mathcal{I})$ | III.293 | relative $\mathcal{I}$-cell complexes |
| $f \square g$ | III.293 | pushout product of $f$ and $g$ |

## Chapter III.8

| | | |
|---|---|---|
| $\mathsf{FinSet}$ | III.300 | category of pointed finite sets and pointed morphisms |
| $\mathcal{F}$ | III.300 | small skeleton of $\mathsf{FinSet}$ with objects $\underline{n} = \{0, 1, \ldots, n\}$ |
| $a^\flat$ | III.301 | punctured finite set |
| $\lvert a \rvert$ | III.301 | cardinality of $a$ |
| $L : \underline{m} \wedge \underline{n} \cong \underline{mn}$ | III.301 | lexicographic order bijection |
| $X : (\mathcal{F}, \underline{0}) \longrightarrow (\mathsf{C}, *)$ | III.302 | $\Gamma$-object in $\mathsf{C}$ |
| $\Gamma$-$\mathsf{C}$ | III.302 | category of $\Gamma$-objects in $\mathsf{C}$ |
| $p_n : X\underline{n} \longrightarrow \prod X\underline{1}$ | III.303 | $n$th Segal map |
| $\Gamma$-$\mathsf{sSet}$ | III.303 | category of $\Gamma$-simplicial sets |
| $\Gamma$-$\mathsf{Cat}$ | III.303 | 2-category of $\Gamma$-categories |
| $N_* : \Gamma\text{-}\mathsf{Cat} \longrightarrow \Gamma\text{-}\mathsf{sSet}$ | III.304 | functor induced by composition with $N$ |
| $\overline{S}$ | III.304 | the $\mathcal{F}$-sphere |
| $h_i, \eta_{m,n,i}, \eta_{m,n}$ | III.304 | structure morphisms for $\underline{m} \wedge X\underline{n} \longrightarrow X\underline{mn}$ |
| $\mathsf{K}^{\mathcal{F}} X$ | III.305 | $K$-theory symmetric spectrum of $\Gamma$-simplicial set $X$ |
| $(\mathsf{C}, \rho) = \{\mathsf{C}_s, \rho_{s,t}\}$ | III.306 | $\underline{n}$-system in $\mathsf{C}$ with $(s,t)$-gluing morphsims $\rho_{s,t}$ |
| $\{\alpha_s\}$ | III.307 | morphism of $\underline{n}$-systems |
| $\mathsf{C}^{\mathcal{F}}\underline{n} = \mathsf{C}^{\mathcal{F}}_{\cong}\underline{n}$, $\mathsf{C}^{\mathcal{F}}_{\mathrm{lax}}\underline{n}$, $\mathsf{C}^{\mathcal{F}}_{\mathrm{co}}\underline{n}$ | III.308 | pointed categories of strong, lax, and colax $\underline{n}$-systems in $\mathsf{C}$ |
| $\mathsf{C}^{\mathcal{F}}\psi$ | III.309 | pointed functor associated to $\psi : \underline{n} \longrightarrow \underline{m}$ |
| $\mathsf{C}^{\mathcal{F}} = \mathsf{C}^{\mathcal{F}}_{\cong}$, $\mathsf{C}^{\mathcal{F}}_{\mathrm{lax}}$, $\mathsf{C}^{\mathcal{F}}_{\mathrm{co}}$ | III.311 | strong, lax, and colax Segal $\Gamma$-categories of $\mathsf{C}$ |

## Chapter III.9

## Chapter III.10

## Chapter III.11

## Chapter III.12

## Chapter III.13

| | | |
|---|---|---|
| $\alpha_F$ | III.494 | canonical map from the homotopy colimit to the colimit |
| $*$ | III.494 | constant functor $\mathsf{Mon}^n(k) \longrightarrow \mathsf{Top}$ |

## Appendix III.A

| | | |
|---|---|---|
| $\mathsf{Bi}^{\mathsf{tsy}}$ | III.511 | full sub-2-category of $\mathsf{Bi}^{\mathsf{sy}}$ with tight objects |
| $\mathsf{SMB}$ | III.511 | tricategory of symmetric monoidal bicategories |
| $\mathsf{SBB}$ | III.512 | tricategory of symmetric bimonoidal bicategories |
| $\mathsf{BB}$ | III.512 | tricategory of bimonoidal bicategories |
| $\mathsf{Bi}^{\mathsf{stsy}}$ | III.522 | strict symmetric bimonoidal categories |
| $(B\Sigma)^+$ | III.522 | Quillen's +-construction of $B\Sigma$ |

# Index

**B**

**C**

**N**

**P**

**U**